\documentclass[12pt]{article}
\title{Poisson deformations of affine symplectic varieties II}
\author{Yoshinori Namikawa}
\date{ }
\usepackage{amsmath,amssymb,amsthm, amscd}
\chardef\bslash=`\\

\def\0{{\mathcal O}}
\def\g{{\mathfrak g}}
\def\q{{\mathfrak q}}
\def\p{{\mathfrak p}}
\def\h{{\mathfrak h}}
\def\k{{\mathfrak k}}

\def\b{{\mathfrak b}}

\begin{document} 
\maketitle 

\begin{center}
{\bf Introduction} 
\end{center}

Let $Y$ be an affine symplectic variety of dimension $2n$, and 
let $\pi: X \to Y$ be a crepant resolution. 
By the definition, there is a symplectic 2-form $\bar{\sigma}$ on 
the smooth part $Y_{reg} \cong \pi^{-1}(Y_{reg})$, and it extends to 
a 2-form $\sigma$ on $X$. Since $\pi$ is crepant, $\sigma$ is a symplectic 
2-form on $X$.  The symplectic structures on $X$ and $Y$ define Poisson structures 
on them in a natural manner. One can define a Poisson deformation of 
$X$ (resp. $Y$) (cf. [Na 1]). A Poisson deformation of $X$ is equivalent 
to a symplectic deformation, namely a deformation of the pair $(X, \sigma)$. 
Let $$\mathrm{PD}_Y: (Art)_{\mathbf C} \to (Set)$$ be the Poisson deformation functor from 
the category of local Artinian {\bf C}-algebras with residue field $\mathbf{C}$ 
to the category of sets (cf. [Na 2], (1.1)). In [Na 2] we have studied a 
morphism of functors $$\pi_*: \mathrm{PD}_X \to \mathrm{PD}_Y$$ induced by 
$\pi$ (cf. [Na 2], \S 5). In particular, $\mathrm{PD}_X$ and $\mathrm{PD}_Y$ are 
both unobstructed and $\pi_*$ is a {\em finite covering}. 
In order to apply these results to geometric situations, 
we need the algebraizations of various formal objects. For this purpose, we start with an affine symplectic 
variety $Y$ with a $\mathbf{C}^*$-actions. More precisely,     
we assume that this ${\mathbf C}^*$-action has a  
unique fixed point $0 \in Y$ and the cotangent space 
$m_{Y,0}/m^2_{Y,0}$ is decomposed into 1-dimensional eigen-spaces with 
{\em positive} weights. Moreover $\bar{\sigma}$ is assumed to be positively 
weighted (cf. [Na 1], (A.1)).   
This $\mathbf{C}^*$-action extends to a 
$\mathbf{C}^*$-action on $X$ ([Na 1], Proposition (A.7), Step 1). By [Na 2]  
one can construct a $\mathbf{C}^*$-equivariant commutative diagram 

\begin{equation} 
\begin{CD} 
\mathcal{X} @>>> \mathcal{Y} \\ 
@VVV @VVV \\ 
\mathbf{A}^d @>>> \mathbf{A}^d     
\end{CD} 
\end{equation} 
where $\mathcal{X} \to \mathbf{A}^d$ (resp. $\mathcal{Y} \to \mathbf{A}^d$) 
is a Poisson deformation of $X = \mathcal{X}_0$ 
(resp. $Y = \mathcal{Y}_0$), and both of them are universal 
at $0 \in \mathbf{A}^d$. Let $\mathrm{PDef}(X)$ (resp. $\mathrm{PDef}(Y)$) be a 
small open neighborhood of $0$ in the 1-st affine space $\mathbf{A}^d$ 
(resp. 2-nd affine space $\mathbf{A}^d$). We call them the Kuranishi 
spaces for the Poisson deformations of $X$ and $Y$, respectively. 
The map $\mathbf{A}^d \to \mathbf{A}^d$ restricts to the map 
$f: \mathrm{PDef}(X) \to \mathrm{PDef}(Y)$. $f$ is a finite surjective map between smooth 
varieties of the same dimension. 

The main result of this paper 
claims that $f$ is a Galois covering (cf. Theorem (1.1)). \S 1 is devoted to 
the proof of Theorem (1.1). For a projective symplectic 
variety $Y$ and its crepant resolution $X$, Markman [Ma] already proved the same result 
for the usual Kuranishi spaces, where he pointed out the Weyl groups of the 
folded Dynkin diagrams appear as the Galois group. A main motivation of this paper was  
to generalize his result for a Poisson deformation of an affine symplectic variety. 
While trying to understand his result, the author realized that his result can be proved 
very naturally in terms of Poisson deformations. This point of view also enables us to reprove 
his original result in a slightly different manner (cf. (1.3)).

In the next section, we apply Theorem (1.1) to the Poisson deformations of an affine 
symplectic variety related to a nilpotent orbit in a complex simple Lie algebra.  
Let $\g$ be a complex simple Lie algebra and let $G$ be the adjoint group. 
For a parabolic subgroup $P$ of $G$, denote by $T^*(G/P)$ the cotangent bundle of 
$G/P$. The image of the Springer map $s: T^*(G/P) \to \g$ is the closure $\bar{O}$ of a nilpotent 
(adjoint) orbit  $O$ in $\g$. Then the normalization $\tilde{O}$ of $\bar{O}$ is an 
affine symplectic variety with the Kostant-Kirillov 2-form. If $s$ is birational onto its image, 
then the Stein factorization $T^*(G/P) \to \tilde{O} \to \bar{O}$ of $s$ gives a crepant 
resolution of $\tilde{O}$. In this situation, we shall prove that the Brieskorn-Slodowy 
diagram (cf. \S 2) 
\begin{equation}
\begin{CD} 
G \times^P r(\p) @>>> \widetilde{G\cdot r(\p)}\\ 
@VVV @VVV \\ 
\k(\p) @>>> \k(\p)/W' 
\end{CD}
\end{equation}  
coincides with the $\mathbf{C}^*$-equivariant commutative 
diagram of the universal Poisson deformations of $T^*(G/P)$ and $\tilde{O}$. 
The precise definitions and notations for the Brieskorn-Slodowy diagram can 
be found in \S 2. Here the group $W'$ appears as the Galois group of the finite 
covering of Kuranishi spaces. The group $W'$ is defined as  
$N_W(L)/W(L)$ (cf. \S 2), where $L$ is a Levi subgroup of $P$, $W(L)$ is the 
Weyl group of $L$, and $N_W(L)$ is the normalizer group of $L$ in the Weyl group 
$W$ of $G$.  This $W'$ coincides with the $W'$ in [Ho]. Howlett has extensively studied 
$W'$ in [Ho]. According to [Ho], $W'$ is {\em almost} a reflection group. 
But, in our situation where the Springer map has degree 1, we can give a geometric proof that 
$W'$ is a reflection group (Lemma (2.2)).  
The author would like to thank E. Markman for pointing out a mistake in 
the proof of (1.1) in the first version. Finally he thanks the referee for reading the 
manuscript very carefully and giving him valuable suggestions.   
\vspace{0.15cm}

{\bf Terminologies}. 
(i) A {\em symplectic variety} $(X, \omega)$ is a pair of a normal algebraic variety $X$ 
defined over {\bf C}  
and a symplectic 2-form $\omega$ on the regular part $X_{reg}$ of $X$ such 
that, for any resolution $\mu: \tilde{X} \to X$, the 2-form $\omega$ on $\mu^{-1}(X_{reg})$ 
extends to a closed regular 2-form on $\tilde{X}$.  We also have a similar notion of a symplectic variety in 
the complex analytic category (eg. the germ of a normal complex space, a holomorphically 
convex, normal, complex space).     
The symplectic 2-form $\omega$ defines a bivector $\Theta \in \wedge^2 \Theta_{X_{reg}}$ 
by the identification $\Omega^2_{X_{reg}} \cong \wedge^2 \Theta_{X_{reg}}$ by $\omega$. 
Define a Poisson structure  
$\{ \;, \; \}$ on $X_{reg}$ by $\{f,g\} := \Theta(df \wedge dg)$. Since $X$ is normal, the  
Poisson structure on $X_{reg}$ uniquely extends to a Poisson structure on $X$.

(ii) Let $X$ be a Poisson variety, and let  $f: \mathcal{X} \to T$ be a 
Poisson deformation of $X$ such that $\mathcal{X}_0 = X$ with $0 \in T$. 
In this paper, we say that $f$ is {\em universal} (or more precisely, {\em formally universal})  
at $0 \in T$ if, for any Poisson deformation 
$\mathcal{X}' \to S$ of $X$ ($\mathcal{X}'_0 = X$, $0 \in S$) with a local Artinian base 
$S$, there is a unique map $(S,0) \to (T,0)$ such that 
$\mathcal{X}' \cong \mathcal{X} \times_T S$ as the Poisson deformations of $X$ over $S$.    
In this case, for a small open neighborhood $V$ of $0 \in T$, the family $f\vert_{f^{-1}(V)}: 
f^{-1}(V) \to V$ is called the Kuranishi family for the Poisson deformations of $X$, and 
$V$ is called the Kuranishi space for the Poisson deformations of $X$. 
\vspace{0.2cm}

\begin{center}
{\bf 1. The Kuranishi spaces for Poisson deformations and Galois coverings}  
\end{center}
\vspace{0.2cm}

Let $(X, \sigma)$ and $(Y, \bar{\sigma})$ be the same as in Introduction. 
Let $\Sigma$ be the singular locus of $Y$. According to [Ka], $\Sigma$ 
is stratified into symplectic varieties. In particular, each stratum has 
even dimension. If $\mathrm{Codim}_Y\Sigma = 2$, then the maximal strata 
parametrize ADE singularities. More precisely, there is a closed 
subset $\Sigma_0 \subset \Sigma$ with $\mathrm{Codim}_Y\Sigma_0 
\geq 4$ and $Y$ is locally isomorphic to 
$(S,0) \times ({\mathbf C}^{2n-2},0)$ at every point $p \in \Sigma - \Sigma_0$,  
where $S$ is an ADE surface singularity. 
Let $\mathcal{B}$ be the set of connected components of 
$\Sigma - \Sigma_0$. Let $B \in \mathcal{B}$. Pick a point $b \in B$ and 
take a transversal slice $S_B \subset Y$ of $B$ passing through $b$.  In other 
words, $Y$ is locally isomorphic to $S_B \times (B, b)$ around $b$. 
$S_B$ is a surface with an ADE singularity. Put $\tilde{S}_B := \pi^{-1}(S_B)$. 
Then $\tilde{S}_B$ is a minimal resolution of $S_B$.  Put $T_B := S_B \times (B,b)$ 
and $\tilde{T}_B := \pi^{-1}(T_B)$. Note that $\tilde{T}_B = \tilde{S}_B \times (B,b)$.  
We put $\bar{\sigma}_B 
:= \bar{\sigma}\vert_{(T_B)_{reg}}$, and $\sigma_B := 
\sigma\vert_{\tilde{T}_B}$. Then $(T_B, \bar{\sigma}_B)$ is a singular 
symplectic variety, and $(\tilde{T}_B, \sigma_B)$ is a smooth symplectic 
variety.  Let $C_i$ ($1 \le i \le r$) be the $(-2)$-curves contained in $\tilde{S}_B$  
and let $[C_i] \in H^2(\tilde{S}_B, \mathbf{R})$ be their classes in the 2-nd 
cohomology group.  Then $$\Phi := \{ \Sigma a_i[C_i]; a_i \in \mathbf{Z}, \; (\Sigma a_i[C_i])^2 = -2\}$$ 
is a root system of the same type as that of the ADE-singularity $S_B$.  Let $W$ 
be the Weyl group of $\Phi$.
 Let $\{E_i(B)\}_{1\le i \le \bar{r}}$ be the set of irreducible exceptional divisors of $\pi$ lying over 
$B$, and let $e_i(B) \in H^2(X, \mathbf{Z})$ be their classes. 
Clearly, $\bar{r} \le r$. If $\bar{r} = r$, then we define $W_B := W$. 
If $\bar{r} < r$, the Dynkin diagram of $\Phi$ has a 
non-trivial graph automorphism. When $\Phi$ is of type $A_r$ with $r > 1$, 
 $\bar{r} = [r+1/2]$ and the Dynkin diagram has a graph automorphism $\tau$ of  order $2$.   
When $\Phi$ is of type $D_r$ with $r \ge 5$, $\bar{r} = r-1$ and the Dynkin diagram 
has a graph automorphism $\tau$ of order $2$.  
When $\Phi$ is of type $D_4$, the Dynkin diagram has two different graph automorphisms 
of order $2$ and $3$. There are two possibilities of $\bar{r}$; $\bar{r} = 2$ or $\bar{r} = 3$. 
In the first case, let $\tau$ be the graph automorphism of order $3$. In the latter case, 
let $\tau$ be the graph automorphism of order $2$. 
Finally, when $\Phi$ is of type $E_6$, $\bar{r} = 4$ and the Dynkin diagram has 
a graph automorphism $\tau$ of order $2$. 
In all these cases, we define 
$$ W_B := \{w \in W; \tau w \tau^{-1} = w \}. $$      

{\bf Theorem (1.1)}. {\em $f: \mathrm{PDef}(X) \to \mathrm{PDef}(Y)$ is a Galois covering with $G = \prod_{B \in \mathcal{B}} W_B$.} 
\vspace{0.15cm}

{\em Proof of Theorem}. We divide the proof in 4 steps. 

{\bf 1} ({\em Outline of the proof}):  In {\bf 2, (ii)} we shall construct the Kuranishi space $\mathrm{PDef}(T_B)$ 
for the Poisson deformation of $T_B$ and the Kuranishi space $\mathrm{PDef}(\tilde{T}_B)$ 
for the Poisson deformation of $\tilde{T}_B$. 
By the construction there is a finite Galois map $f_B: \mathrm{PDef}(\tilde{T}_B) \to 
\mathrm{PDef}(T_B)$, whose Galois group 
is the Weyl group of the root system corresponding to the ADE-singularity $S_B$. 
Since $T_B$ (resp. $\tilde{T}_B$) is an open set of $Y$ (resp. $X$), any 
Poisson deformation of $Y$ (resp. $X$) over a local Artinian base induces a 
Poisson deformation of $T_B$ (resp. $\tilde{T}_B$) over the same base. Thus we have 
a morphism of functors $\mathrm{PD}_Y \to \mathrm{PD}_{T_B}$ 
(resp. $\mathrm{PD}_X \to \mathrm{PD}_{\tilde{T}_B}$). 
Since $R^1\pi_*\mathcal{O}_X = 0$ and $\pi_*\mathcal{O}_X = \mathcal{O}_Y$, 
the crepant resolution $\pi: X \to Y$ induces a morphism of functors 
$\mathrm{PD}_X \to \mathrm{PD}_Y$ (cf. [Na 2], Proof (i) of Theorem (5.1)). 
Similarly, we have a morphism of functors 
$\mathrm{PD}_{\tilde{T}_B} \to \mathrm{PD}_{T_B}$. 
These morphisms form a commutative diagram  

\begin{equation} 
\begin{CD} 
\mathrm{PD}_X @>>> \mathrm{PD}_{\tilde{T}_B} \\ 
@VVV @VVV \\ 
\mathrm{PD}_Y @>>> \mathrm{PD}_{T_B}.     
\end{CD} 
\end{equation} 

For a complex space $W$ with an origin $0 \in W$, we denote 
by $\hat{W}$ the formal completion of $W$ at $0$. 
The commutative diagram above induces a commutative diagram of 
formal spaces

\begin{equation} 
\begin{CD} 
\widehat{\mathrm{PDef}(X)} @>>> \widehat{\mathrm{PDef}(\tilde{T}_B)} \\ 
@VVV @VVV \\ 
\widehat{\mathrm{PDef}(Y)} @>>>  \widehat{\mathrm{PDef}(T_B)}.      
\end{CD} 
\end{equation} 

By using the period maps (cf. {\bf 2} (i)), one can see that
this diagram is actually the formal completions of a commutative diagram of 
complex spaces (cf. {\bf 4} (i))

\begin{equation} 
\begin{CD} 
\mathrm{PDef}(X) @>{\varphi_B}>> \mathrm{PDef}(\tilde{T}_B) \\ 
@VVV @V{f_B}VV \\ 
\mathrm{PDef}(Y) @>>>  \mathrm{PDef}(T_B).      
\end{CD} 
\end{equation} 
  
Put $V_B := \mathrm{Im}(\varphi_B)$. We shall 
prove that (a) $V_B$ and $f_B(V_B)$ are both non-singular, and 
(b) $f_B\vert_{V_B}: V_B \to f_B(V_B)$ is a $W_B$-Galois covering (cf. {\bf 3}, 
and the final part of {\bf 4}, (ii)).  
Then we get the following commutative diagram 

\begin{equation} 
\begin{CD} 
\mathrm{PDef}(X) @>{\varphi_B}>> \prod_{B \in \mathcal{B}} V_B \\ 
@VVV @VVV \\ 
\mathrm{PDef}(Y) @>>> \prod_{B \in \mathcal{B}} f_B(V_B).      
\end{CD} 
\end{equation} 

We finally prove that the induced map 
$$\iota: \mathrm{PDef}(X) \to \mathrm{PDef}(Y) \times_{\prod f_B(V_B)} 
\prod V_B$$ is an isomorphism (cf. {\bf 4}, (iii)). 

{\bf 2} ({\em Poisson deformations and period maps}):  
(i) Let us consider the commutative diagram of universal Poisson deformations in the 
introduction

\begin{equation} 
\begin{CD} 
\mathcal{X} @>>> \mathcal{Y} \\ 
@V{\alpha}VV @VVV \\ 
\mathbf{A}^d @>>> \mathbf{A}^d.     
\end{CD} 
\end{equation} 

Note that $\mathcal{Y}$ has a $\mathbf{C}^*$-action with positive weights with a 
fixed point $0 \in Y$. Moreover, the diagram is $\mathbf{C}^*$-equivariant and 
$\alpha: \mathcal{X} \to \mathbf{A}^d$ 
is a simultaneous resolution of $\mathcal{Y} \to \mathbf{A}^d$. Then we see that 
$\mathcal{X}$ is a $C^{\infty}$-trivial fiber bundle over $\mathbf{A}^d$  
by [Slo 2], Remark at the end of section 4.2\footnote{In [Slo 2] 4.2, $\mathcal{Y}$ is assumed to 
be smooth. But the arguments can be applied to our case.} . 
Let  $\Omega^{\cdot}_{\mathcal{X}^{an}/\mathbf{A}^d}$ be the relative complex-analytic de Rham 
complex. Let $\mathcal{K}$ be the subsheaf of $\Omega^2_{\mathcal{X}^{an}/\mathbf{A}^d}$ 
which consists of d-closed relative 2-forms. By the natural map $\mathcal{K}[-2] \to 
\Omega^{\cdot}_{\mathcal{X}^{an}/\mathbf{A}^d}$, we can define a sequence of maps: 
$$ \alpha_*\mathcal{K} \to  \mathbf{R}^2\alpha_*\Omega^{\cdot}_{\mathcal{X}^{an}/\mathbf{A}^d} 
\cong R^2\alpha_*\alpha^{-1}\mathcal{O}^{an}_{\mathbf{A}^d}.$$ 
Since $R^2\alpha_*\alpha^{-1}\mathcal{O}^{an}_{\mathbf{A}^d} \cong R^2\alpha_*\mathbf{C} \otimes_{\mathbf C}
\mathcal{O}^{an}_{\mathbf{A}^d}$ (cf. [Lo], Lemma (8.2)), we have an isomorphism 
$$R^2\alpha_*\alpha^{-1}\mathcal{O}^{an}_{\mathbf{A}^d} \cong H^2(X, \mathbf{C}) \otimes_{\mathbf C} 
\mathcal{O}^{an}_{\mathbf{A}^d}.$$ Since $\alpha$ is a Poisson deformation of $X$, 
$\mathcal{X}$ admits a relative symplectic 2-form $\sigma_{\mathcal X}$ such that 
$\sigma_{\mathcal X}\vert_X = \sigma$. Then $\sigma_{\mathcal X}$ gives a section $s$ of the 
sheaf $H^2(X, \mathbf{C}) \otimes_{\mathbf C} 
\mathcal{O}^{an}_{\mathbf{A}^d}$.  Let $$ev_t: H^2(X, \mathbf{C}) \otimes_{\mathbf C} 
\mathcal{O}^{an}_{\mathbf{A}^d} \to H^2(X, \mathbf{C})$$ be the evaluation map at 
$t \in \mathbf{A}^d$. We define a period map 
$$ p: \mathbf{A}^d \to H^2(X, \mathbf{C}) $$ by $p(t) = ev_t(s)$. By the construction, 
$p$ is a holomorphic map. In [G-K], Proposition (5.4), one can find another approach to the definition of 
the period map. 
The period map restricts to give a map $$p_X: \mathrm{PDef}(X)  
\to H^2(X, \mathbf{C}).$$ We also call this map the period map for $\alpha$. 
Since the tangential map $T_0 \mathrm{PDef}(X) \to H^2(X, \mathbf{C})$ 
of $p_X$ is an isomorphism ([Na 1], Corollary 10), $p_X$ is an open immersion.
 
(ii) Since $S_B$ is an ADE-singularity, it is 
isomorphic to the hypersurface defined by a weighted homogeneous polynomial $f(x,y,z)$. 
Then $S_B$ has a ${\mathbf C}^*$-action with positive weights, and 
$\bar{\tau}_B:= Res(dx \wedge dy \wedge dz/f)$ is a generator of $K_{S_B}$ with a positive weight. 
By the minimal resolution $\tilde{S}_B \to S_B$, $\bar{\tau}_B$ is pulled back to a 
symplectic 2-form $\tau_B$ on $\tilde{S}_B$. On the other hand, $(B,b)$ is isomorphic to 
$(\mathbf{C}^{2n-2},0)$ and it admits a canonical symplectic 2-form $\tau_{\mathbf{C}^{2n-2}} 
:= ds_1 \wedge dt_1 + ... + ds_{n-1} \wedge dt_{n-1}$ with the standard coordinates 
$(s_1, ..., s_{n-1}, t_1, ..., t_{n-1})$ 
of $\mathbf{C}^{2n-2}$. 
By a generalization of Darboux's theorem ([Na 2], Lemma 1.3), one can see that 
$(T_B, \bar{\sigma}_B)$ is equivalent to $(T_B, \bar{\tau}_B + \tau_{\mathbf{C}^{2n-2}})$ as a symplectic variety. 
Therefore, $(\tilde{T}_B, \sigma_B)$ is equivalent to  
$(\tilde{T}_B, {\tau}_B + \tau_{\mathbf{C}^{2n-2}})$. Let $\g$ be the complex simple Lie algebra 
of the same type as $S_B$ and let $\h$ be a Cartan subalgebra of $\g$. We denote by 
$W$ the Weyl group of $\g$. By using a special transversal slice of $\g$, one can 
construct the universal Poisson deformation $$\mathcal{S}_B \to \h/W$$ of $S_B$ and 
the universal Poisson deformation $$\tilde{\mathcal{S}}_B \to \h$$ of $\tilde{S}_B$ ([Na 2], Proposition (3.1), (1)).  
Moreover, the composite $${\mathcal S}_B \times (\mathbf{C}^{2n-2},0) \stackrel{p_1}\to 
{\mathcal S}_B \to \h/W $$ and the composite $$\tilde{\mathcal S}_B \times (\mathbf{C}^{2n-2},0) 
\stackrel{p_1}\to \tilde{\mathcal S}_B \to \h$$ are respectively universal Poisson deformations 
of $T_B$ and $\tilde{T}_B$ (ibid, Proposition (3.1), (2)).      
Note that $T_B$ has a $\mathbf{C}^*$-action with positive weights. Since $\tilde{S}_B$ is 
the minimal resolution of $S_B$, this $\mathbf{C}^*$-action uniquely extends to that on 
$\tilde{T}_B$. We now have a $\mathbf{C}^*$-equivariant commutative diagram  

\begin{equation} 
\begin{CD} 
\tilde{\mathcal S}_B \times (\mathbf{C}^{2n-2},0) @>>> {\mathcal S}_B \times (\mathbf{C}^{2n-2},0) \\ 
@VVV @VVV \\ 
\h @>>> \h/W.     
\end{CD} 
\end{equation} 

As is seen in (i), one can define a period map 
$\h \to H^2(\tilde{T}_B, \mathbf{C})$.  
The Kuranishi space $\mathrm{PDef}(\tilde{T}_B)$ for the Poisson deformation of  
$\tilde{T}_B$ is an open neighborhood of $0 \in \h$ and the period map above 
restricts  to a map $$p_B: \mathrm{PDef}(\tilde{T}_B) \to H^2(\tilde{T}_B, \mathbf{C}) 
\cong H^2(\tilde{S}_B, \mathbf{C}).$$  
Since the tangential map $T_0 \mathrm{PDef}(\tilde{T}_B) \to H^2(\tilde{T}_B, \mathbf{C})$ 
of $p_B$ is an isomorphism ([Na 1], Corollary 10), $p_B$ is an open immersion.  
Since $\mathrm{PDef}(X)$ (resp. $\mathrm{PDef}(\tilde{T}_B)$) can be regarded 
as open subsets of $H^2(X, \mathbf{C})$ (resp. $H^2(\tilde{T}_B, \mathbf{C})$) 
by the period maps, one can define a holomorphic map 
$$\varphi_B: \mathrm{PDef}(X) \to \mathrm{PDef}(\tilde{T}_B) $$ 
so that the following diagram commutes 

\begin{equation} 
\begin{CD} 
\mathrm{PDef}(X) @>{p_X}>> H^2(X, \mathbf{C}) \\ 
@V{\varphi_B}VV @VVV \\ 
\mathrm{PDef}(\tilde{T}_B)   @>{p_B}>> H^2(\tilde{S}_B, \mathbf{C}).      
\end{CD} 
\end{equation} 

Since a Poisson deformation of $X$ over a local Artinian base restricts 
to give a Poisson deformation of $\tilde{T}_B$ over the same base, we 
have a natural morphism of functors $\mathrm{PD}_X \to \mathrm{PD}_{\tilde{T}_B}$. 
By the (formal) universality of $\mathrm{PDef}(\tilde{T}_B)$ a formal map 
$$\hat{\varphi}_B: \widehat{\mathrm{PDef}(X)} \to \widehat{\mathrm{PDef}(\tilde{T}_B)}$$ 
is uniquely determined. This is nothing but the formal completion of $\varphi_B$.  
 
{\bf 3} ({\em Description of $\mathrm{Im}(\varphi_B)$}): 
The Weyl group $W$ (of the root system $\Phi$ associated with $\tilde{S}_B$) acts on 
$H^2(\tilde{S}_B, \mathbf{C})$.
The period map $p: \h \to H^2(\tilde{S}_B, \mathbf{C})$ is a 
$W$-equivariant linear map by [Ya] (cf. [Na 2] Proof of Proposition (3.2)). 
The commutative diagram  

\begin{equation} 
\begin{CD} 
\h @>{p}>> H^2(\tilde{S}_B, \mathbf{C}) \\ 
@VVV @VVV \\ 
\h/W @>>> H^2(\tilde{S}_B, \mathbf{C})/W      
\end{CD} 
\end{equation} 
induces a commutative diagram 

\begin{equation} 
\begin{CD} 
\mathrm{PDef}(\tilde{T}_B) @>{p_B}>> H^2(\tilde{S}_B, \mathbf{C}) \\ 
@V{f_B}VV @V{q}VV \\ 
\mathrm{PDef}(T_B) @>>> H^2(\tilde{S}_B, \mathbf{C})/W.      
\end{CD} 
\end{equation} 
 
Let $E_i(B)$ ($1 \le i \le \bar{r}$) be the irreducible exceptional divisors of $\pi: X \to Y$ lying over $B$, 
and let $e_i(B) \in H^2(X, \mathbf{Z})$ be the cohomology class determined by $E_i(B)$. 
Even if $E_i(B)$ is irreducible, $E_i(B) \cap \tilde{S}_B$ might be reducible.   
Denote by $r_B$ the restriction map $H^2(X, \mathbf{Z}) \to H^2(\tilde{S}_B, 
\mathbf{Z})$. Then $\mathrm{Im}(r_B \otimes \mathbf{C})$ is the $\mathbf{C}$-vector 
subspace generated by $r_B(e_i(B))$, ($1 \le i \le \bar{r}$).  
Let $r$ be the number of  $(-2)$-curves in $\tilde{S_B}$.  As explained at the beginning 
of this section, when $\bar{r}  <  r$, 
we have a non-trivial graph automorphism $\tau$  
of the Dynkin diagram.  We then define $\Gamma_B := <\tau>$.  
When $\bar{r} = r$, we just put $\Gamma_B = id$. Then  
$H^2(\tilde{S}_B, \mathbf{C})^{\Gamma_B} = \mathrm{Im}(r_B \otimes \mathbf{C})$.  
We put $\tilde{V}_B := \mathrm{Im}(r_B \otimes \mathbf{C})$. 
Let $W'$ be the subgroup of $W$ consisting of the elements which preserve 
$\tilde{V}_B$ as a set, and let $W_B :=\{w \in W; \tau w {\tau}^{-1} = w\}$. 
This $W_B$ is nothing but $W^1$ in 
[Ca], Chapter 13. 
It is obvious that $W_B \subset W'$.   
\vspace{0.15cm}

{\bf Lemma (1.2)}. $W_B = W'$. 
\vspace{0.15cm}

{\em Proof}. 
For simplicity we put $V := H^2(\tilde{S}_B, \mathbf{C})$ and 
$V^{\tau} := H^2(\tilde{S}_B, \mathbf{C})^{\Gamma_B}$. 
Let $(V^{\tau})^{\perp}$ be the orthogonal complement of $V^{\tau}$ with 
respect to the inner product. Assume that $\tau^2 = 1$. Assume that 
$g \in W$ preserves $V^{\tau}$ as a set. Since $g$ is an isometry of $V$, it 
acts on $(V^{\tau})^{\perp}$. Since $\tau$ acts on 
$(V^{\tau})^{\perp}$ by $-1$, we see that $g$ commutes with $\tau$. 
This means that $g \in W_B$. 
We next treat the case where $\tau$ has order $3$ (ie. $D_4$ case). 
By the first part in [Ca], Chapter 13, $W$ is normalized by 
$\tau$; in other words, $\tau W \tau^{-1} = W$ in the 
automorphism group of the root system. Assume that $W_B$ does not 
coincide with $W'$. Then there is an element $w' \in W'$ such that $w:= \tau w' {\tau}^{-1}$ 
does not equal $w'$. Note that  $w' w^{-1} \ne 1$ and $w' w^{-1}$ acts trivially 
on $V^{\tau}$. We shall show that such an element does not exist. 
Let $V = (\mathbf{C}^4, (\;, \;))$ be a 4-dimensional real vector space with a positive definite 
symmetric form. Let $e_1, ..., e_4$ be an orthogonal basis such that $(e_i, e_j) = 0$ if 
$i \ne j$ and $(e_i,e_i) = 1$. One can choose simple roots for $D_4$ in such a way 
that  $C_1 := e_1 - e_2$, $C_2 := e_2 - e_3$, $C_3 := e_3 - e_4$ and 
$C_4 := e_3 + e_4$.  Define $\tau$ by $\tau(C_1) = c_3$, $\tau(C_3) = C_4$ 
and $\tau(C_4) = C_1$. Then $V^{\tau}$ is a 2-dimensional vector space spanned 
by $e_2 - e_3$ and $e_1-e_2 + 2e_3$. 
Every element $w$ of the Weyl group has the form 
$$ w(e_i) = (-1)^{\epsilon_i}e_{\sigma(i)} $$ with a permutation $\sigma$ of 
$1$, $2$, $3$ and $4$. Here each $\epsilon_i$ is $0$ or $1$ and 
$\Sigma \epsilon_i$ is even. It is easily checked that if $w$ acts on $V^{\tau}$ 
trivially, then $w = id$.   Q.E.D.         
\vspace{0.15cm}

The natural map $\tilde{V}_B \to q(\tilde{V}_B)$ factors through $\tilde{V}_B/W_B$.  
By Lemma (1.2) the map $\tilde{V}_B/W_B \to q(\tilde{V}_B)$ is the normalization 
map.  Since $W_B$ is generated by reflections 
(cf. [Ca], Chapter 13), $\tilde{V}_B/W_B$ is smooth.  
Put $V_B := \tilde{V}_B \cap \mathrm{PDef}(\tilde{T}_B)$. 
Then there is a commutative diagram          

\begin{equation} 
\begin{CD} 
V_B @>>> \mathrm{PDef}(\tilde{T}_B)\\ 
@VVV @V{f_B}VV \\ 
f_B(V_B) @>>> \mathrm{PDef}(T_B)      
\end{CD} 
\end{equation} 

We shall prove that the image of the map $\mathrm{PDef}(X) \to \mathrm{PDef}(\tilde{T}_B)$  
coincides with $V_B$. Since 
$\tilde{V}_B = \mathrm{Im}(r_B \otimes \mathbf{C})$, the image is contained in 
$V_B$ by the commutative diagram of the period 
maps $p_X$ and $p_B$. The tangent space $T_{([X,\sigma])} \mathrm{PDef}(X)$ (resp. 
$T_{([\tilde{T}_B, \sigma_B])} \mathrm{PDef}(\tilde{T}_B)$) is identified with 
$H^2(X, \mathbf{C})$ (resp. $H^2(\tilde{S}_B, \mathbf{C})$). Moreover,  
the tangential map $$T_{([X,\sigma])} \mathrm{PDef}(X) \to    
T_{([\tilde{T}_B, \sigma_B])} \mathrm{PDef}(\tilde{T}_B)$$ is identified with 
the map $$ r_B \otimes \mathbf{C}: 
H^2(X, \mathbf{C}) \to H^2(\tilde{S}_B, \mathbf{C}).$$ 
This means that 
the image of the map $\mathrm{PDef}(X) \to \mathrm{PDef}(\tilde{T}_B)$  
coincides with $V_B$. 
 
{\bf 4} ({\em Proof that $\iota$ is an isomorphism}): 
(i) There is a commutative diagram of functors 

\begin{equation} 
\begin{CD} 
\mathrm{PD}_X @>>> \mathrm{PD}_{\tilde{T}_B} \\ 
@VVV @VVV \\ 
\mathrm{PD}_Y @>>> \mathrm{PD}_{T_B}.      
\end{CD} 
\end{equation} 

Corresondingly we have a commutative 
diagram of formal spaces 

\begin{equation} 
\begin{CD} 
\widehat{\mathrm{PDef}(X)} @>{\hat{\varphi}}_B>> \widehat{\mathrm{PDef}(\tilde{T}_B)} \\ 
@V{\hat{f}}VV @V{\hat{f}_B}VV \\ 
\widehat{\mathrm{PDef}(Y)} @>{\hat{\phi}_B}>> \widehat{\mathrm{PDef}(T_B)}.     
\end{CD} 
\end{equation} 

We have seen in {\bf 2} that the formal maps $\hat{f}$, $\hat{f}_B$ and $\hat{\varphi_B}$ 
are completions of the holomorphic maps $f$, $f_B$ and $\varphi_B$. 
Let us prove that $\hat{\phi}_B$ is also the completion of a holomorphic map 
$\phi_B: \mathrm{PDef}(Y) \to \mathrm{PDef}(T_B)$. 
In fact, there is a commutative diagram of local rings 

\begin{equation} 
\begin{CD} 
\mathcal{O}_{\mathrm{PDef}(T_B),0} @>{\hat{\phi}_B^*\vert_{\mathrm{PDef}(T_B),0}}>> \hat{\mathcal{O}}_{\mathrm{PDef}(Y),0} \\ 
@VVV @VVV \\ 
\mathcal{O}_{\mathrm{PDef}(\tilde{T}_B),0} @>{\varphi_B^*}>> \hat{\mathcal{O}}_{\mathrm{PDef}(X),0}      
\end{CD} 
\end{equation} 

Since $\mathrm{Im}(\varphi_B^*) \subset \mathcal{O}_{\mathrm{PDef}(X),0}$, 
we see that 
$$\hat{\phi}_B^*(\mathcal{O}_{\mathrm{PDef}(T_B),0}) \subset 
\mathcal{O}_{\mathrm{PDef}(X),0} \cap \hat{\mathcal{O}}_{\mathrm{PDef}(Y),0}.$$ 
On the right hand side, we take the intersection in 
$\hat{\mathcal{O}}_{\mathrm{PDef}(X),0}$. 
We shall prove that $$\mathcal{O}_{\mathrm{PDef}(X),0} \cap \hat{\mathcal{O}}_{\mathrm{PDef}(Y),0} 
= \mathcal{O}_{\mathrm{PDef}(Y),0}.$$ 
For simplicity, we put $A = \mathcal{O}_{\mathrm{PDef}(Y),0}$ and $B = \mathcal{O}_{\mathrm{PDef}(X),0}$. 
Let $m$ be the maximal ideal of $A$. Note that $B$ is a finite $A$-module. 
Assume that $g \in \hat{A} \cap B$. For $n >0$, we can write $g = g_n + h_n$ with 
$g_n \in A$ and $h_n \in m^n\hat{A}$. Since $h_n = g - g_n \in B$, we have 
$$h_n \in m^n\hat{A} \cap B \subset m^n\hat{B} \cap B = m^nB.$$
In other words, $g \in A + m^nB$. Since $n$ is arbitrary, we have 
$$g \in \cap_{n > 0}(A + m^nB) = A.$$ This shows that $\hat{\phi}_B$ induces 
a homomorphism $\mathcal{O}_{\mathrm{PDef}(T_B),0} \to 
\mathcal{O}_{\mathrm{PDef}(Y),0}$, hence a holomorphic map  
$$\phi_B: \mathrm{PDef}(Y) \to \mathrm{PDef}(T_B).$$ 
As a consequence we have a commutative diagram 

\begin{equation} 
\begin{CD} 
\mathrm{PDef}(X) @>{\varphi_B}>> \mathrm{PDef}(\tilde{T}_B) \\ 
@V{f}VV @V{f_B}VV \\ 
\mathrm{PDef}(Y) @>{\phi_B}>> \mathrm{PDef}(T_B).     
\end{CD} 
\end{equation} 

(ii) We shall briefly recall some results proved in [Na 2]. 
Put $U := Y \setminus \Sigma_0$ and $\tilde{U} := \pi^{-1}(U)$. 
There are natural morphisms of functors $\mathrm{PD}_Y \to \mathrm{PD}_U$ 
and $\mathrm{PD}_X \to \mathrm{PD}_{\tilde U}$. By Lemma (5.3) of [Na 2] these 
are both isomorphisms. Denote by $\mathbf{PT}^1_U$ (resp. $\mathbf{PT}^1_{\tilde U}$) 
the tangent space of $\mathrm{PD}_U$ (resp. $\mathrm{PD}_{\tilde U}$). 
We have an isomorphism $\mathbf{PT}^1_{\tilde U} \cong 
H^2(\tilde{U}, \mathbf{C})$ (cf. Proof (i) of [Na 2], Theorem (5.1)). 
In [Na 2] we have constructed a local system $\mathcal{H}$ of $\mathbf{C}$-modules 
on $\Sigma - \Sigma_0$. The local system $\mathcal{H}$ is the subsheaf of 
$\underline{\mathrm{Ext}}^1(\Omega^1_U, \mathcal{O}_U)$ which consists of local 
sections coming from the Poisson deformations ([Na 2], (1.4) and (1.5)). 
We have an exact sequence (cf. [Na 2], Proposition (1.11)) 
$$0 \to H^2(U, \mathbf{C}) \to \mathbf{PT}^1_U \to H^0(\Sigma - \Sigma_0, \mathcal{H}).$$ 
Here the first term $H^2(U, \mathbf{C})$ is the space of locally trivial Poisson 
deformations of $U$.   
There is a commutative diagram of exact sequences 

\begin{equation} 
\begin{CD} 
0 @>>> H^2(U, \mathbf{C}) @>>> H^2(\tilde{U}, \mathbf{C}) @>>> 
H^0(U, R^2(\pi_{\tilde U})_*\mathbf{C}) \\ 
@. @VVV @VVV @.\\ 
0 @>>> H^2(U, \mathbf{C}) @>>> \mathbf{PT}^1_U @>>> H^0(\Sigma - \Sigma_0, \mathcal{H}).      
\end{CD} 
\end{equation} 
 
Theorem (5.1) of [Na 2] claims that $\mathrm{PD}_U$ and $\mathrm{PD}_{\tilde U}$ 
are both unbostructed and have the same dimension. In the course of its proof, 
we also prove that two maps $$H^2(\tilde{U}, \mathbf{C}) \to  
H^0(U, R^2(\pi_{\tilde U})_*\mathbf{C})$$ and 
$$\mathbf{PT}^1_U \to H^0(\Sigma - \Sigma_0, \mathcal{H})$$ are both surjective. 
Suppose that $\pi\vert_{\tilde U}: \tilde{U} \to U$ has exactly $m$ irreducible 
components. Then we have  
$$\dim \mathrm{Im}[H^2(\tilde{U}, \mathbf{C}) \to  
H^0(U, R^2(\pi_{\tilde U})_*{\mathbf C})] = m$$ and   
$$\dim \mathrm{Im}[\mathbf{PT}^1_U \to H^0(\Sigma - \Sigma_0, \mathcal{H})] = m.$$
There is a natural injection  
$H^0(U, R^2(\pi_{\tilde U})_*{\mathbf C}) 
\to \prod H^2(\tilde{S}_B, \mathbf{C})$ and the image of the composed map 
$H^2(\tilde{U}, \mathbf{C}) \to  
H^0(U, R^2(\pi_{\tilde U})_*{\mathbf C}) \to 
\prod H^2(\tilde{S}_B, \mathbf{C})$ has dimension $m$.
Since $H^2(X, \mathbf{C}) \cong H^2(\tilde{U}, \mathbf{C})$ (cf.[Na 2], Lemma (5.3)), 
$$ \dim \mathrm{Im}[H^2(X, \mathbf{C}) \to \prod H^2(\tilde{S}_B, \mathbf{C})] = m.$$ 
On the other hand, the left hand side equals $\Sigma \dim V_B$ by the argument in {\bf 3}. 
Hence we have $m = \Sigma \dim V_B$. 
Denote by  $T^1_{S_B}$ 
the tangent space of $Def(S_B)$ at the origin. 
The stalk $\mathcal{H}_b$ of $\mathcal{H}$ at $b \in B$ is 
isomorphic to $T^1_{S_B}$ (cf. [Na 2], (1.3), (1.5), (3.1)).  
Since $\mathcal{H}$ is a local system, there is 
a natural injection 
$H^0(\Sigma - \Sigma_0, \mathcal{H}) \to \prod_{B \in \mathcal{B}}T^1_{S_B}$. 
The image of the composed map $\mathbf{PT}^1_U \to H^0(\Sigma - \Sigma_0, \mathcal{H}) \to 
\prod_{B \in \mathcal{B}}T^1_{S_B}$ has dimension $m$.
Note that $T^1_{S_B}$ is identified with the tangent space $T_0\mathrm{PDef}(T_B)$  
of $\mathrm{PDef}(T_B)$ ([Na 2], Proposition (3.1)).   
The map $\mathbf{PT}^1_U \to \prod T^1_{S_B}$ coincides with 
$T_0 \mathrm{PDef}(Y) \to \prod T_0 \mathrm{PDef}(T_B)$.
Hence we have 
$$ \dim \mathrm{Im}[T_0\mathrm{PDef}(Y) \to \prod T_0\mathrm{PDef}(T_B)] = m.$$ 
On the other hand, the image of the map $$\mathrm{PDef}(Y) \to \prod \mathrm{PDef}(T_B)$$  
is $\prod f_B(V_B)$. 
Since $\mathrm{PDef}(Y)$ is smooth and $m = \Sigma \dim f_B(V_B)$, we see that $f_B(V_B)$ is smooth. 
In particular, the map 
$V_B \to  f_B(V_B)$ is a finite Galois cover with Galois group $W_B$. 
\vspace{0.2cm}

(iii) Let us consider the commutative diagram 

\begin{equation} 
\begin{CD} 
\mathrm{PDef}(X) @>{\prod \varphi_B}>> \prod_{B \in \mathcal{B}} \mathrm{PDef}(\tilde{T}_B) \\ 
@V{f}VV @V{\prod f_B}VV \\ 
\mathrm{PDef}(Y) @>{\prod \phi_B}>> \prod_{B \in \mathcal{B}} \mathrm{PDef}(T_B).     
\end{CD} 
\end{equation} 

By the argument in {\bf 3}, the horizontal map at the bottom factorizes as  
$\mathrm{PDef}(Y) \to \prod_{B \in \mathcal{B}}f_B(V_B) \to 
\prod \mathrm{PDef}(T_B)$ and 
the following diagram commutes:

\begin{equation} 
\begin{CD} 
\mathrm{PDef}(X) @>>> \prod_{B \in \mathcal{B}} V_B \\ 
@VVV @VVV \\ 
\mathrm{PDef}(Y) @>>> \prod_{B \in \mathcal{B}} f_B(V_B)      
\end{CD} 
\end{equation} 
 
The commutative diagram above induces a map 
$$\iota: \mathrm{PDef}(X) \to \mathrm{PDef}(Y) \times_{\prod f_B(V_B)} 
\prod V_B.$$ 
First we shall prove that $\mathrm{PDef}(Y) \times_{\prod f_B(V_B)} 
\prod V_B$ is smooth. 
Since 
$m = \Sigma \dim f_B(V_B)$ and each $f_B(V_B)$ is smooth, the map 
$\mathrm{PDef}(Y) \to \prod f_B(V_B)$ is a smooth map. Therefore, 
$\mathrm{PDef}(Y) \times_{\prod f_B(V_B)} \prod V_B \to \prod V_B$ is also 
a smooth map, and $\mathrm{PDef}(Y) \times_{\prod f_B(V_B)} \prod V_B$ is 
smooth.        
Finally we shall prove that the map $\iota$ is an isomorphism. 
The tangent space $T$ of $\mathrm{PDef}(Y) \times_{\prod f_B(V_B)} \prod V_B$ at the 
origin $\{0\} \times \prod \{0\}$ is isomorphic to 
$\mathbf{PT}^1_U \times_{\prod T_0(f_B(V_B))} \prod T_0 V_B$. 
Since $T_0(V_B) \to T_0(f_B(V_B))$ is the zero map, it is isomorphic to   
$H^2(U, \mathbf{C}) \oplus \prod V_B$. The map $d\iota: T_0 \mathrm{PDef}(X)(\cong H^2(\tilde{U}, \mathbf{C}))  
\to T$ is injective. 
In fact, if $v \in H^2(\tilde{U}, \mathbf{C})$ is mapped to zero by this map, 
then $v$ must be sent to zero by the map $H^2(\tilde{U}, \mathbf{C}) 
\to T_0 V_B$ for each $B$. In other words, $v$ is sent to zero by the 
composed map $H^2(\tilde{U}, \mathbf{C}) \to H^0(U, R^2(\pi_{\tilde U})_*{\mathbf C}) 
\to \prod H^2(\tilde{S}_B, \mathbf{C})$. Since the second map 
$H^0(U, R^2(\pi_{\tilde U})_*{\mathbf C}) 
\to \prod H^2(\tilde{S}_B, \mathbf{C})$ is an injection, 
$v$ is already sent to zero by the first map 
$H^2(\tilde{U}, \mathbf{C}) \to H^0(U, R^2(\pi_{\tilde U})_*{\mathbf C})$. 
By the commutative diagram above, we see that
$v \in H^2(U, \mathbf{C})$. On the other hand, $v$ must be sent to zero 
by the map $H^2(\tilde{U}, \mathbf{C}) \to \mathbf{PT}^1_U$. The restriction of 
this map to $H^2(U, \mathbf{C})$ is an injection; hence $v = 0$. 
Since both $T_0 \mathrm{PDef}(X)$ and $T$ have the same dimension 
$h^2(U, \mathbf{C}) + m$, 
$d\iota$ is an isomorphism.   
Note that $\mathrm{PDef}(X)$ and 
$\mathrm{PDef}(Y) \times_{\prod f_B(V_B)} \prod V_B$ are 
both smooth; hence $\iota$ is an isomorphism.  (End of the proof of Theorem (1.1)) 
\vspace{0.2cm} 

{\bf Remark}.  
For $g \in W_B$, denote by $g_B : \mathrm{PDef}(\tilde{T}_B) \to \mathrm{PDef}(\tilde{T}_B)$ the automorphism induced by $g$. 
By Theorem (1.1), $g$ also induces an automorphism $g_X: \mathrm{PDef}(X) \to \mathrm{PDef}(X)$. 
By pulling back the universal family $\tilde{\mathcal{T}}_B \to \mathrm{PDef}(\tilde{T}_B)$ by $g_B : 
\mathrm{PDef}(\tilde{T}_B) \to \mathrm{PDef}(\tilde{T}_B)$, we have a new family 
$\tilde{\mathcal{T}}'_B \to \mathrm{PDef}(\tilde{T}_B)$. 
Since $G$ acts trivially on the universal family $\mathcal{T}_B \to \mathrm{PDef}(T_B)$, 
we have a diagram of the birational maps 
$$ \tilde{\mathcal{T}}_B \rightarrow \mathcal{T}_B \times_{\mathrm{PDef}(T_B)} \mathrm{PDef}(\tilde{T}_B) 
\leftarrow \tilde{\mathcal{T}}'_B. $$ If $g \ne 1$, the birational map 
$\tilde{\mathcal{T}}_B --\to \tilde{\mathcal{T}}'_B$ is not regular. 
Similarly, by pulling back the universal family $\mathcal{X} \to \mathrm{PDef}(X)$ by $g_X$, we 
get a diagram of the birational maps 
$$ \mathcal{X} \rightarrow \mathcal{Y} \times_{\mathrm{PDef}(Y)} \mathrm{PDef}(X) \leftarrow \mathcal{X}'.$$ 
If $g \ne 1$, the birational map $\mathcal{X} --\to \mathcal{X}'$ is not regular; moreover, 
some $\pi$-exceptional divisor $E \subset X$ lying over $B$ is contained in its indeterminacy locus.   
\vspace{0.2cm}

{\bf  (1.3)}. The proof of Theorem (1.1) also works for the compact case. 
More exactly, let $Y$ be a projective symplectic variety and let $\pi: X \to Y$ be 
a crepant (projective) resolution. We assume that $h^0(X, \Omega^2_X) = 1$ and 
$h^1(X, \mathcal{O}_X) = 0$. 
By [Na], the Kuranishi spaces $\mathrm{Def}(X)$ and $\mathrm{Def}(Y)$ are both non-singular, and 
the induced map $\bar{f}: \mathrm{Def}(X) \to \mathrm{Def}(Y)$ is a finite surjective map.  
Markman [Ma] proved that $\bar{f}$ is actually a Galois covering whose Galois 
group coincides with $\prod W_B$. 
Let $\mathrm{PDef}(X)$ be the Kuranishi space for the Poisson deformations
(symplectic deformations) of $(X, \sigma)$. 
Let $\mu: \mathcal{X} \to \mathrm{Def}(X)$ be the Kuranishi family in the usual sense. 
Then $\mathcal{V} := \mu_*\Omega^2_{\mathcal{X}/\mathrm{Def}(X)}$ is a line bundle on $\mathrm{Def}(X)$, 
and $\mathcal{V}^* := \mathcal{V} - \{0-section\}$ is a ${\mathbf C}^*$-bundle over $\mathrm{Def}(X)$. 
The Kuranishi space $\mathrm{PDef}(X)$ for the Poisson deformation of $(X, \sigma)$ is defined as an open 
neighborhood of $(X, \sigma) \in \mathcal{V}^*$. In particular, $\mathrm{PDef}(X)$ is a smooth variety.  
One can define a period map $$ p_X:  \mathrm{PDef}(X) \to H^2(X, \mathbf{C})$$ by $p_X(X_t, \sigma_t) 
:= [\sigma_t] \in H^2(X, \mathbf{C})$.  Let $Q \subset H^2(X, {\mathbf C})$ be the hypersurface 
defined by  $q = 0$ with the Beauville-Bogomolov form $q$ ([Be]). 
Let $\bar{Q} \subset \mathbf{P}(H^2(X, \mathbf{C}))$ be the projective hypersurface 
defined by $q = 0$. There is a commutative diagram 

\begin{equation} 
\begin{CD} 
\mathrm{PDef}(X) @>{p_X}>> Q -\{0\} \\ 
@VVV @VVV \\ 
\mathrm{Def}(X) @>{\bar{p}_X}>> \bar{Q}      
\end{CD} 
\end{equation} 
where $\bar{p}_X$ is the usual period map (cf. [Be]). The fibers of 
both vertical maps are $\mathbf{C}^*$ and $p_X$ maps the fibers isomorphically. 
Since $\bar{p}_X$ is an open immersion by the local Torelli theorem, 
$p_X$ is also an open immersion by the commutative diagram.    
As in (1.1), {\bf 2}, we have a 
commutative diagram of period maps  

\begin{equation} 
\begin{CD} 
\mathrm{PDef}(X) @>>> \mathrm{PDef}(\tilde{T}_B) \\ 
@VVV @VVV \\ 
H^2(X, \mathbf{C}) @>>> H^2(\tilde{S}_B, \mathbf{C})      
\end{CD} 
\end{equation} 

Note that $[\sigma] \in H^2(X, \mathbf{C})$ is not zero, but 
$[\sigma_B] = 0$ in $H^2(\tilde{S}_B, \mathbf{C})$(cf. [Ka], Corollary 2.8) .    
Let $\mathrm{PD}_Y$ be the Poisson deformation functor of $(Y, \bar{\sigma})$ (cf. [Na 2]).  
As in [Na 2], we shall prove that $\mathrm{PD}_Y$ is unobstructed. 
Let $U := Y - \Sigma_0$ and put $\tilde{U}:= \pi^{-1}(U)$.
The key commutative diagram in the compact case is    

\begin{equation} 
\begin{CD} 
0 @>>> \mathbf{H}^2(U, \tilde{\Omega}^{\geq 1}_U) @>>> \mathbf{H}^2(\tilde{U}, \Omega^{\geq 1}_{\tilde{U}}) @>>> 
H^0(U, R^2(\pi_{\tilde U})_*\mathbf{C}) \\ 
@. @VVV @VVV @.\\ 
0 @>>> \mathbf{H}^2(U, \tilde{\Omega}^{\geq 1}_U) @>>> \mathbf{PT}^1_U @>>> H^0(\Sigma - \Sigma_0, \mathcal{H}).      
\end{CD} 
\end{equation} 

The exact sequence on the first row comes from the Leray spectral sequence 
$$ E^{p,q}_2 := H^p(U, R^q(\pi_{\tilde U})_*\Omega^{\geq 1}_{\tilde U}) \Rightarrow 
\mathbf{H}^{p+q}(\tilde{U}, \Omega^{\geq 1}_{\tilde U}). $$ 
In the exact sequence on the second row, $\mathbf{H}^2(U, \tilde{\Omega}^{\geq 1}_U)$ 
is the space of  the 1-st order Poisson deformations of $U$ which are locally trivial as 
flat deformations (cf. the proof of [Na 2], Lemma (1.9)), and $\mathbf{PT}^1_U$ 
is the tangent space of $\mathrm{PD}_U$. 
Note that $\mathbf{H}^2(\tilde{U}, \Omega^{\geq 1}_{\tilde{U}})$ is the tangent space of 
$\mathrm{PD}_{\tilde U}$ (cf. [Na 1], Propositions 8, 9).   
Since $H^2(X, \mathbf{C}) \cong H^2(\tilde{U}, \mathbf{C})$ (cf. the proof of Proposition 2 
of [Na 3]), $H^2(\tilde{U})$ has a pure Hodge structure of weight 2.  By the distinguished triangle 
$$ \Omega^{\geq 1}_{\tilde U} \to \Omega^{\cdot}_{\tilde U} \to \mathcal{O}_{\tilde U} 
\stackrel{[1]} \to \Omega^{\geq 1}_{\tilde U}[1], $$ and by the fact that $H^1(\tilde{U}, \mathcal{O}_{\tilde U}) 
= 0$,  we have an exact sequence 
$$ 0 \to \mathbf{H}^2(\tilde{U}, \Omega^{\geq 1}_{\tilde{U}}) \to H^2(\tilde{U}, \mathbf{C}) 
\to H^2(\tilde{U}, \mathcal{O}_{\tilde U}). $$  
By the proof of Proposition 2 of [Na 3], $Gr^0_F(H^2({\tilde U})) = H^2(\tilde{U}, \mathcal{O}_{\tilde U})$. 
Therefore, $\mathbf{H}^2(\tilde{U}, \Omega^{\geq 1}_{\tilde{U}}) = F^1(H^2({\tilde U}))$.   
Since $F^1(H^2(X)) = H^{2,0}(X) \oplus H^{1,1}(X)$, we see that 
the map $$\mathbf{H}^2(\tilde{U}, \Omega^{\geq 1}_{\tilde{U}}) \to 
H^0(U, R^2(\pi_{\tilde U})_*\mathbf{C})$$ is identified with the map 
$$ H^{2,0}(X) \oplus H^{1,1}(X) \subset H^2(X, \mathbf{C}) \to 
H^0(U, R^2(\pi_{\tilde U})_*\mathbf{C}).$$ 
Since $U$ has only quotient singularities, 
$U$ is {\bf Q}-factorial. Then, by [Ko-Mo], Proposition (12.1.6), 
$$\mathrm{Im}[H^2(X, \mathbf{C}) \to H^0(U, R^2(\pi\vert_{\tilde U})_*\mathbf{C})] 
= \mathrm{Im}[\Sigma_{B,i} \mathbf{C}[E_i(B)] \to H^0(U, R^2(\pi\vert_{\tilde U})_*\mathbf{C})].$$ 
Since the map $H^0(U, R^2(\pi\vert_{\tilde U})_*\mathbf{C}) \to \prod H^2(\tilde{S}_B, \mathbf{C})$ is an injection, 
we see that $$\dim \mathrm{Im}[H^2(X, \mathbf{C}) \to H^0(U, R^2(\pi\vert_{\tilde U})_*\mathbf{C})] 
= \dim \mathrm{Im}[H^2(X, \mathbf{C}) \to \prod H^2(\tilde{S}_B, \mathbf{C})] = m, $$ where 
$m = \Sigma \dim V_B$. 
By the map $H^2(X, \mathbf{C}) \to H^0(Y, R^2\pi_*\mathbf{C})$, $H^{0,2}(X)$ is sent to zero.
Therefore $$\dim \mathrm{Im}[\mathbf{H}^2(\tilde{U}, \Omega^{\geq 1}_{\tilde{U}}) \to 
H^0(U, R^2(\pi_{\tilde U})_*\mathbf{C})] = m.$$ 
By [Na 2], Proposition (4.2), $\dim H^0(\Sigma - \Sigma_0, \mathcal{H}) = m$. 
Then, the same argument in [Na], Theorem (4.1) can be applied to our case 
to prove that $\mathrm{PD}_Y$ is unobstructed. 
Since $H^0(Y, \Theta_Y) = 0$, $\mathrm{PD}_Y$ is pro-represented by a complete regular local 
ring $R$ over $\mathbf{C}$. We denote by $\widehat{\mathrm{PDef}(Y)}$
the formal 
scheme \footnote{The author has not yet constructed the Kuranishi space for $\mathrm{PD}_Y$ 
as a complex space. Let $\beta : \mathcal{Y} \to \mathrm{Def}(Y)$ be the 
universal family and let $\mathcal{Y}^0 \subset \mathcal{Y}$ be the locus where $\beta$ is 
smooth. We put $\beta^0 := \beta\vert_{\mathcal{Y}^0}$. Then 
$\bar{\mathcal V} := (\beta^0)_*\Omega^2_{\mathcal{Y}^0/\mathrm{Def}(Y)}$ 
seems very likely to be a line bundle on ${\mathrm{Def}(Y)}$  
(cf. [Na 6]).  Then the Kuranishi space for $\mathrm{PD}_Y$ would be realized as an open subset of 
$\bar{\mathcal V}-\{0-section\}$} defined by $R$. Since $R^1\pi_*\mathcal{O}_X = 0$ and $\pi_*\mathcal{O}_X 
= \mathcal{O}_Y$, the crepant resolution $\pi : X \to Y$ induces a morphism 
of functors $\mathrm{PD}_X \to \mathrm{PD}_Y$ (cf. [Na 2], Proof (i) of Theorem (5.1)). 
By the formal universality of $\widehat{\mathrm{PDef}(Y)}$, we have a formal map 
$\widehat{\mathrm{PDef}(X)} \to \widehat{\mathrm{PDef}(Y)}$. Note that 
$\dim \widehat{\mathrm{PDef}(X)} = \dim \widehat{\mathrm{PDef}(Y)}$, and  
as in [Na], Lemma (4.2), the formal map is finite.  There is a commutative diagram 
\begin{equation} 
\begin{CD} 
\widehat{\mathrm{PDef}(X)} @>>> \widehat{\mathrm{PDef}(Y)} \\ 
@VVV @VVV \\ 
\widehat{\mathrm{Def}(X)} @>>> \widehat{\mathrm{Def}(Y)}      
\end{CD} 
\end{equation} 
 
The fibers of the maps $\widehat{\mathrm{PDef}(X)} \to \widehat{\mathrm{Def}(X)}$ 
and $\widehat{\mathrm{PDef}(Y)} \to \widehat{\mathrm{Def}(Y)}$ both have dimension 
1, and they correspond to the Poisson deformations of $X$ and $Y$ where the underlying 
flat deformations are fixed and the symplectic structures only vary.  
Therefore, 
$$ \widehat{\mathrm{PDef}(X)} \cong \widehat{\mathrm{PDef}(Y)} \times_{\widehat{\mathrm{Def}(Y)}} 
\widehat{\mathrm{Def}(X)}.$$  
In order to prove Markman's result, we only have to prove the 
map $\widehat{\mathrm{PDef}(X)} \to \widehat{\mathrm{PDef}(Y)}$ is a finite Galois cover with 
Galois group $\prod_{B \in \mathcal{B}}W_B$.  
The rest of the argument is similar to (1.1).  
Another approach avoiding formal schemes is the following. 
We first remark that Poisson deformations of $\tilde{T}_B$ and $T_B$ are 
consequently determined only by the underlying flat deformations (cf. {\bf 2}, (ii) and 
[Na 2], Proposition (3.1)). 
In particular, the map $\varphi_B : \mathrm{PDef}(X) \to \mathrm{PDef}(\tilde{T}_B)$ 
factorizes as $$ \mathrm{PDef}(X) \to \mathrm{Def}(X) \to \mathrm{PDef}(\tilde{T}_B).$$ 
By the same reasoning as (1.1), we then have a commutative diagram 

\begin{equation} 
\begin{CD} 
\mathrm{Def}(X) @>>> \prod \mathrm{PDef}(\tilde{T}_B) \\ 
@VVV @VVV \\ 
\mathrm{Def}(Y) @>>> \prod \mathrm{PDef}(T_B).       
\end{CD} 
\end{equation} 

The induced map $$\mathrm{Def}(X) \cong \mathrm{Def}(Y) \times_{\prod \mathrm{PDef}(T_B)} 
\prod \mathrm{PDef}(\tilde{T}_B)$$ turns out to be an isomorphism. 
\vspace{0.2cm}

\begin{center}
{\bf 2. Poisson deformations associated with nilpotent orbits} 
\end{center}

Let $\g$ be a complex simple Lie algebra.  We fix a Cartan 
subalgebra $\h$ of $\g$, and let $\p$ be a parabolic subalgebra 
of $\g$ such that $\h \subset \p$.  Denote by $r(\p)$ (resp. $n(\p)$) 
the solvable radical (resp. nilpotent radical) of $\p$. Define 
$\k(\p) := \h \cap r(\p)$, and let $l(\p)$ be the Levi subalgebra of 
$\p$ which contains $\h$. Let $G$ be the adjoint group of $\g$, and 
$P$ (resp. $L$) the closed subgroup of $G$ corresponding 
to $\p$ (resp. $l(\p)$). The cotangent bundle $T^*(G/P)$ of $G/P$ is 
isomorphic to the vector bundle $G \times^P n(\p)$ over $G/P$. 
The Springer map $s: G \times^P n(\p) \to \g$ is defined by 
$s([g,x]) := Ad_g(x)$ for $[g,x] \in G \times^P n(\p)$. 
The image $\mathrm{Im}(s)$ is the closure $\bar{O}$ of a nilpotent orbit $O$. 
The Springer map is a generically finite, surjective, projective morphism. 
Let $G \times^P n(\p) \to \tilde{O} \to \bar{O}$ be the Stein factorization of $s$.  
Let $W$ be the Weyl group of $G$. If we fix a Borel subalgebra $\b$ such that 
$\h \subset \b \subset \p$, then $\p$ is determined by a choice of a subset $J$ of 
the set of simple roots ([Na 4], (P1)). Then $W(L)$ is generated by reflections in 
elements of $J$; hence $W(L)$ is a subgroup of $W$. Define  
$W' := N_W(L)/W(L)$, where $N_W(L)$ is the normalizer group of $L$ 
in $W$, and $W(L)$ is the Weyl group of $L$. $N_W(L)$ acts on $\k(\p)$, where 
$W(L)$ acts trivially on $\k(\p)$. Therefore, $W'$ acts effectively on $\k(\p)$.   
Let us construct the Brieskorn-Slodowy diagram. 
There is 
a direct sum decomposition 
$$ r(\p) = \mathfrak{k}(\p) 
\oplus n(\p), (x \to x_1 + x_2)$$ 
where $n(\p)$ is the 
nil-radical of $\mathfrak{p}$.
We have a well-defined map 
$$ G \times^P r(\p) \to 
\mathfrak{k}(\p) $$ 
by sending 
$[g, x] \in G \times^P r(\p)$ to 
$x_1 \in \mathfrak{k}(\p)$
([Slo], 4.3). 
On the other hand, define a map 
$G \times^P r(\p) \to G\cdot r(\p)$ 
by $[g,x] \to Ad_g(x)$. By the adjoint quotient map 
$\g \to \h/W$, we have a map from $G\cdot r(\p)$ to 
$\h/W$. These maps form a commutative diagram ([Slo], 4.3) 
 
\begin{equation}
\begin{CD} 
G \times^P r(\p) @>>> G\cdot r(\p)\\ 
@V{f}VV @VVV \\ 
\k(\p) @>>> \h/W. 
\end{CD}
\end{equation}  
One can find an instructive example of the diagram in [Na 5], Example 7.10.
 
Let $\widetilde{G\cdot r(\p)}$ be the normalization of 
$G\cdot r(\p)$. Here the set $G\cdot r(\p)$ is irreducible since 
it is the image of the smooth variety $G \times^P r(\p)$;  we regard 
$G\cdot r(\p)$ as a variety 
with the reduced structure. The normalization of the image of the map 
$G \cdot r(\p) \to \h/W$ coincides with $\k(\p)/W'$. 
Then the map $\widetilde{G\cdot r(\p)} \to \h/W$ factors through 
$\k(\p)/W'$.   
By [Na 4], Lemma (1.1),  we already know that 
$G\cdot r(\p) \times_{\h/W} \k(\p)$ is irreducible.  
Let $\widetilde{G\cdot r(\p) \times_{\h/W} \k(\p)}$ 
be the normalization of the variety  
$(G\cdot r(\p) \times_{\h/W} \k(\p))_{red}$.
Then $\widetilde{G\cdot r(\p)}$ is the quotient variety of 
$\widetilde{G\cdot r(\p) \times_{\h/W} \k(\p)}$ by $W'$. 
The variety $\widetilde{G\cdot r(\p) \times_{\h/W} \k(\p)}$ has 
a resolution $G \times^P r(\p)$ whose canonical line bundle is 
trivial  ([Na 4], Lemma 1.2). In particular, 
$\widetilde{G\cdot r(\p) \times_{\h/W} \k(\p)}$ has only rational 
singularities (cf. Ibid, Lemma 1.2). 
Hence its quotient variety $\widetilde{G\cdot r(\p)}$ 
also has rational singularities. In particular, $\widetilde{G\cdot r(\p)}$ is 
Cohen-Macaulay. 
\vspace{0.2cm}

{\bf Lemma (2.1)}. 
{\em  
The central fiber $F$ of 
$\widetilde{G\cdot r(\p) \times_{\h/W} \k(\p)} \to \k(\p)$ is isomorphic 
to $\tilde{O}$.}

{\em Proof}. Since $\widetilde{G\cdot r(\p) \times_{\h/W} \k(\p)}$ 
is Cohen-Macaulay and $\k(\p)$ is smooth, the central fiber $F$ is 
also Cohen-Macaulay. On the other hand, let us consider the birational map 
$G \times^P r(\p) \to \widetilde{G\cdot r(\p) \times_{\h/W} \k(\p)}$, and 
take their central fibers to get a map 
$T^*(G/P) \to F_{red}$ with connected fibers. Since 
the Springer map is generically finite, this map is birational by Zariski's Main Theorem. 
Moreover, it is an isomorphism outside 
certain codimension 2 subset $Z$ of $F_{red}$.  Take a point $x \in F_{red} - Z$. 
Then we have a surjection 
$$\mathcal{O}_{F,x} \to \mathcal{O}_{F_{red}, x} \cong \mathcal{O}_{T^*(G/P), x}.$$  
By the lemma of Nakayama, this implies that 
$\mathcal{O}_{\widetilde{G\cdot r(\p) \times_{\h/W} \k(\p)}, x} \cong  
\mathcal{O}_{G \times^P r(\p), x}$. 
Therefore, $F$ is reduced at $x$, and moreover, $F$ is smooth at $x$.   
Since $F$ is Cohen-Macaulay and regular in codimension one, $F$ is normal. 
This means that $F = \tilde{O}$. Q.E.D. 
\vspace{0.15cm}

In the remainder, we always assume  
\vspace{0.15cm}

{\bf Assumption}: {\em The Springer map $s: T^*(G/P) \to \bar{O}$ is birational.}
\vspace{0.15cm}

We often write $X_{\p}$ for $G \times^P r(\p)$, and   
$Y_{\k(\p)}$ for $\widetilde{G\cdot r(\p) \times_{\h/W} \k(\p)}$. 
Denote by $\mu_{\p}$ the birational map from $X_{\p}$ to $Y_{\k(\p)}$. 
The map $\mu_{\p}$ is a crepant resolution of $Y_{\k(\p)}$, which is 
an isomorphism in codimension one ([Na 4], Theorem 1.3).  Let $X_{\p,0}$ (resp. $Y_{\k(\p),0}$) 
be the central fiber of $X_{\p} \to \k(\p)$ (resp. $Y_{\k(\p)} \to \k(\p)$). 
Note that $X_{\p,0} = T^*(G/P)$, and $Y_{\k(\p),0} = \tilde{O}$ by Lemma (2.1). 
The birational map $\mu_{\p,0}: X_{\p,0} \to 
Y_{\k(\p),0}$ coincides with the 
Stein factorization of the Springer map $s: T^*(G/P) \to \bar{O}$.  
 
We shall briefly review (P2), (P3), \S 2 of [Na 4]. Let $\mathcal{S}(l(\p))$ be the set of parabolic 
subalgebras $\p'$ which contain $l(\p)$ as Levi subalgebras. 
Then every crepant resolution of $Y_{\k(\p)}$ is isomorphic to $\mu_{\p'}: X_{\p'} 
\to Y_{\k(\p)}$ with $\p' \in \mathcal{S}(l(\p))$ ([Na 4], Theorem (1.3)). Let $M(L) := 
\mathrm{Hom}_{alg.gp.}(L, \mathbf{C}^*)$. 
The 2-nd cohomology group $H^2(X_{\p'}, \mathbf{R})$ is isomorphic to 
$M(L) \otimes \mathbf{R}$. By this isomorphism, 
the nef cone $\overline{\mathrm{Amp}}(\mu_{\p'})$ is regarded 
as a cone in $M(L) \otimes \mathbf{R}$. The cohomology group 
$H^2(X_{\p',0}, \mathbf{R})$ is also isomorphic to 
$M(L) \otimes \mathbf{R}$.  By this isomorphism 
the nef cone $\overline{\mathrm{Amp}}(\mu_{\p',0})$ is regarded 
as a cone in $M(L) \otimes \mathbf{R}$. Note that 
$\overline{\mathrm{Amp}}(\mu_{\p'}) =  \overline{\mathrm{Amp}}(\mu_{\p',0})$. 
One has the following ([Na 4], Remark (1.6)):  
$$M(L) \otimes \mathbf{R} = \cup_{\p' \in \mathcal{S}(l(\p))}\overline{\mathrm{Amp}}(\mu_{\p'}).$$
We say that two nef cones $\overline{\mathrm{Amp}}(\mu_{\p'})$ and 
$\overline{\mathrm{Amp}}(\mu_{\p''})$ are adjacent to each other if they share a common codimension one face. 
In this case, we also say that $\p'$ and $\p''$ are adjacent to each other.  If $\p'$ and $\p''$ are 
adjacent to each other, $\p'$ and $\p''$ are related by an operation called the {\em twist}. 
There are two kinds of twists; a twist of the 1-st kind, and a twist of the 2-nd kind. 
If $\p'$ and $\p''$ are adjacent to each other, $X_{\p'}$ and $X_{\p''}$ are connected by a 
flop (cf. [Na 4], \S 1). If the corresponding twist is of the 1-st kind (resp. 2-nd kind), we say that the flop 
is of the 1-st kind (resp. 2-nd kind).  Let $\mathcal{S}^1(l(\p))$ be the set of parabolic subalgebras in 
$\mathcal{S}(l(\p))$ which can be obtained from $\p$ by a finite succession of the twists of the 
1-st kind. Then $\cup_{\p' \in \mathcal{S}^1(l(\p))}\overline{\mathrm{Amp}}(\mu_{\p'})$ 
coincides with the movable cone $\overline{\mathrm{Mov}}(\mu_{\p,0})$. 
This movable cone $\overline{\mathrm{Mov}}(\mu_{\p,0})$ is a fundamental domain of 
$M(L) \otimes \mathbf{R}$ by the action of $W'$. In particular, 
$$ M(L) \otimes \mathbf{R} = \cup_{w \in W'} w(\overline{\mathrm{Mov}}(\mu_{\p,0})).$$ 
Moreover,  $w(\overline{\mathrm{Mov}}(\mu_{\p,0}))$ coincides with the movable cone 
$\overline{\mathrm{Mov}}(\mu_{w(\p),0})$. One has 
$$ \overline{\mathrm{Mov}}(\mu_{w(\p),0}) = \cup_{\p' \in \mathcal{S}^1(l(\p))} 
\overline{\mathrm{Amp}}(\mu_{w(\p'),0}).$$  
\vspace{0.2cm}

By [Ho] $W'$ is almost a reflection group. But $W'$ turns out to be a reflection group 
under Assumption. \vspace{0.2cm}

{\bf Lemma (2.2)}. {\em The group $W'$ is generated by reflections of $\k(\p)$. In particular, 
$\k(p)/W'$ is smooth.}

{\em Proof}.  Assume that $\p$ and $\p'$ are adjacent to each other, and 
they are related by a 2-nd twist. Then $\p' = w(\p)$ for some $w \in W'$. 
 Let $\phi_w: X_{\p} \to X_{w(\p)}$ be the isomorphism 
defined by $[g,x] \to [gw^{-1}, Ad_w(x)]$ for $[g,x] \in G \times^{P} r(\p)$. 
Let $\bar{\phi}_w : Y_{\k(\p)} \to Y_{\k(\p)}$ be the automorphism induced by 
the map $id \times w : G\cdot r(\p) \times_{\h/W} \k(\p) \to 
G\cdot r(\p) \times_{\h/W} \k(\p)$. Then we have a commutative diagram 
\begin{equation} 
\begin{CD} 
X_{\p} @>{\phi_w}>> X_{w(\p)} \\ 
@VVV @VVV \\ 
Y_{\k(\p)} @>{\bar{\phi}_w}>> Y_{\k(\p)}      
\end{CD} 
\end{equation} 
  
The composite $X_{\p} \to X_{w(\p)} --\to X_{\p}$ induces an automorphism of 
$H^2(X_{\p}, \mathbf{R})$.  We call this automorphism $\varphi_w$.  
In this way, $W'$ acts on $H^2(X_{\p}, \mathbf{R})$. 
Note that its dual action coincides with the natural action of 
$W'$ on $\k(\p)$ by [Na 4], Lemma 2.1.  
We shall prove that $\varphi_w$ is a reflection, that is, 
it is an involution which fixes all points in certain hyperplane.  
The flop $X_{w(\p)} --\to X_{\p}$ can be expressed more exactly by the diagram 
$$ X_{w(\p)} \to 
\widetilde{G \times^{\bar{P}}{\bar{P}}\cdot r(\p)} 
\times_{\k(l(\bar{\p}) \cap \p)/W''}\k(l(\bar{\p}) \cap \p)
\longleftarrow X_{\p}.$$
Here $W''$ is the subgroup of the Weyl group $W(l(\bar{\p}))$ which stabilizes 
$\k(l(\bar{\p}) \cap \p)$ as a set. The element $w$ is contained in $W''$; hence it acts 
on \linebreak  
$\widetilde{G \times^{\bar{P}}{\bar{P}}\cdot r(\p)} 
\times_{\k(l(\bar{\p}) \cap \p)/W''}\k(l(\bar{\p}) \cap \p)$ by $id \times w$. 
Put $$Z_{\bar{\p}} := \widetilde{G \times^{\bar{P}}{\bar{P}}\cdot r(\p)}
\times_{\k(l(\bar{\p}) \cap \p)/W''}\k(l(\bar{\p}) \cap \p)$$ for short.
 
Then we have a commutative diagram  
$$\hspace{2.3cm} X_{\p} \stackrel{\phi_w}\to X_{w(\p)} --\to X_{\p}$$ 
$$ \downarrow \hspace{1.0cm} \downarrow $$ 
$$Z_{\bar{\p}} \stackrel{id \times w}\to Z_{\bar{\p}} $$ 
$$ \downarrow \hspace{1.0cm} \downarrow $$ 
$$Y_{\k(\p)} \to Y_{\k(\p)}.$$   
The automorphism $Z_{\bar{\p}} \stackrel{id \times w}\to Z_{\bar{\p}} $ induces the identity map on 
$H^2(Z_{\bar{\p}}, \mathbf{R})$. 
Since the image of  the map $H^2(Z_{\bar{\p}}, \mathbf{R}) \to H^2(X_{\p}, \mathbf{R})$ has 
codimension one, the automorphism $\varphi_w$ of $H^2(X_{\p}, \mathbf{R})$ is a reflection.  
By [Na 4], Proposition 2.3, $$H^2(X_{\p}, \mathbf{R}) = \bigcup_{w \in W'} 
w(\overline{\mathrm{Mov}}(\mu_{\p,0})).$$ Note that $w(\overline{\mathrm{Mov}}(\mu_{\p,0})) 
= \overline{\mathrm{Mov}}(\mu_{w(\p),0})$. If $w(\overline{\mathrm{Mov}}(\mu_{\p,0}))$ 
and $w'(\overline{\mathrm{Mov}}(\mu_{\p,0}))$ are adjacent to each other, then $w^{-1}w'$ is a reflection. 
For any $w \in W'$, one can connect $w(\overline{\mathrm{Mov}}(\mu_{\p,0}))$ and 
$\overline{\mathrm{Mov}}(\mu_{\p,0})$ by a finite sequence of movable cones 
$w_1(\overline{\mathrm{Mov}}(\mu_{\p,0}))$, $w_2(\overline{\mathrm{Mov}}(\mu_{\p,0}))$, 
..., $w_n(\overline{\mathrm{Mov}}(\mu_{\p,0}))$ with $w_1 = 1$ and $w_n = w$ in such 
a way that $w_i(\overline{\mathrm{Mov}}(\mu_{\p,0}))$ and 
$w_{i+1}(\overline{\mathrm{Mov}}(\mu_{\p,0}))$ are adjacent for all $i$. Then, $w$ 
can be represented as a product of reflections: 
$w = (w_n\cdot w^{-1}_{n-1}) \cdot \cdot\cdot (w_2\cdot w^{-1}_1)$.     
\vspace{0.2cm}
 
{\bf Remark}. When the Springer map $s: T^*(G/P) \to \bar{O}$ 
is not birational, Lemma (2.2) does not hold. See [Na 4], Example 1.9, Remark 2.4. 
\vspace{0.2cm}

{\bf Corollary (2.3)}. {\em $\widetilde{G\cdot r(\p)}$ 
is flat over $\k(\p)/W'$.}

{\em Proof}. First note that $\widetilde{G\cdot r(\p)}$ is Cohen-Macaulay. 
Since every fiber of the map $\widetilde{G\cdot r(\p)} \to \k(\p)/W'$ 
has the dimension equal to $\dim \widetilde{G\cdot r(\p)} - \dim \k(\p)/W'$ and 
$\k(\p)/W'$ is smooth, the map $\widetilde{G\cdot r(\p)} \to \k(\p)/W'$ is flat.   
\vspace{0.2cm}

{\bf Lemma (2.4)}. 
{\em $\widetilde{G\cdot r(\p)} \times_{\k(\p)/W'} \k(\p)$ is a variety. } 
\vspace{0.2cm}

{\em Proof}. 
Let $B$ be 
the affine ring of $\widetilde{G\cdot r(\p)}$, and 
let $A$ (resp. $A'$) be the affine ring of $\k(\p)/W'$ (resp. $\k(\p)$).
Denote by $L$ the quotient field of $B$ and by  
$K$ (resp. $K'$) the quotient field of $A$ (resp. $A'$).  Since 
$B$ is flat over $A$, 
$B \otimes_A A'  \subset B \otimes_A K' = B \otimes_A K \otimes_K K'$. 
Since $B \otimes_A K$ is the localization of $B$ by the multiplicative 
set $S:= A - \{0\}$, it is naturally contained in $L$. Therefore, 
$B \otimes_A K \otimes_K K' \subset L \otimes_K K'$. Since $K'$ is  
a separable extension of $K$, $L \otimes_K K'$ is reduced.  
Since $G\cdot r(\p) \times_{\h/W} \k(\p)$ is irreducible, 
we see that $L \otimes_K K'$ is an integral domain. 
Finally we conclude that $B \otimes_A A'$ is an integral domain 
because  $B \otimes_A A'  \subset L \otimes_K K'$. 
\vspace{0.2cm}

{\bf Lemma (2.5)}. 
{\em  
$$ \widetilde{G\cdot r(\p)} \times_{\k(\p)/W'} \k(\p) \cong 
\widetilde{G\cdot r(\p) \times_{\h/W} \k(\p)}.$$ 
The map  
$\widetilde{G\cdot r(\p)} \to \k(\p)/W'$ is flat, and its 
central fiber coincides with $\tilde{O}$.} 
\vspace{0.2cm}
 
{\em Proof}.   
Let $A$, $A'$ and $B$ be the same as in the proof of Lemma (2.4).  
Consider the map $G\cdot r(\p) \to \h/W$. If $\bar{B}$ and $C$ are 
affine rings of $G\cdot r(\p)$ and $\h/W$ respectively, then 
$\bar{B}$ is a $C$-algebra. The origin $0 \in \h/W$ corresponds to a  
maximal ideal $\bar{m}$ of $C$.  
Let $B'$ be the affine ring 
of  $\widetilde{G\cdot r(\p) \times_{\h/W} \k(\p)}$. By  
Lemma (2.4), we already know that $B \otimes_A A'$ is an integral domain. We have an injection $B \otimes_A A' 
\to B'$. Let $m'$ (resp. $m$) be the maximal ideal of $A'$ (resp. $A$) corresponding to the origin 
$0 \in \k(\p)$ (resp. $0 \in \k(\p)/W'$). 
By the base change property, $B \otimes_A (A'/m') = B /mB$. Then the injection above induces a homomorphism 
$B/mB \to B'/m'B'$. Moreover, there is a map $\bar{B}/\bar{m}\bar{B} \to B/mB$. 
By the definition, $\mathrm{Spec}(\bar{B}/\bar{m}\bar{B})_{red} = \bar{O}$. 
In our case, $\mathrm{Spec}(B'/m'B') = \tilde{O}$ is the normalization of $\bar{O}$. 
Note that the normalization map $\tilde{O} \to \bar{O}$ is an isomorphism in codimension one. 
Therefore, the cokernel of the map $\bar{B}/\bar{m}\bar{B} \to B'/m'B'$ has the support with codimension $\geq 2$ 
in $\tilde{O}$.  The cokernel $Q$ of the map $B/mB \to B'/m'B'$ also has the support with codimension $\geq 2$ 
in $\tilde{O}$.  Take a point $\q \in \tilde{O} - \mathrm{Supp}(Q)$. Then, by the lemma of 
Nakayama, we have an isomorphism  $(B \otimes_A A')_{\q} \cong B'_{\q}$. In particular, 
$\mathrm{Spec}(B/mB)$ and $\mathrm{Spec}(B'/mB')$ are isomorphic in codimension one; hence 
we see that $\mathrm{Spec}(B/mB)$ is regular in codimension one. On the other hand, by Corollary (2.3), $B$ is 
Cohen-Macaulay,  and is flat over $A$.  
Therefore, $\mathrm{Spec}(B/mB)$ is Cohen-Macaulay. 
This means that $\mathrm{Spec}(B/mB)$ is normal; hence $\mathrm{Spec}(B/mB) = \tilde{O}$. 
\vspace{0.2cm}

{\bf Remark}. When the Springer map $s: T^*(G/P) \to \bar{O}$ 
is not birational,  the central fiber of the map  
$G\cdot r(\p) \to \h/W$ is everywhere non-reduced.  
\vspace{0.2cm}

{\bf Proposition (2.6)}. 
{\em Two flat morphisms $$G \times^P r(\p) \to \k(\p)$$ and 
$$\widetilde{G\cdot r(\p)} \to \k(\p)/W'$$ are respectively Poisson deformations  
of $T^*(G/P)$ and $\tilde{O}$.} 
\vspace{0.2cm}

{\em Proof}. The smooth variety $G \times^P r(\p)$ over $\k(\p)$ admits a $G$-invariant 
relative symplectic 2-form $\omega$ ([C-G], Proposition (1.4.14))\footnote{In [C-G], coajoint orbits 
of $\g^*$ are treated. But the Killing form of $\g$ identifies the coadjoint orbits with 
adjoint orbits. The variety $G \times_P (\lambda + \p^{\perp})$ is identified 
with a fiber of $G \times^P r(\p) \to \k(\p)$ by the Killing form.}. 
Let $\omega_t$ be the restriction of $\omega$ to the fiber $G \times^P(t + n(\p))$ over 
$t \in \k(\p)$. There is a $G$-equivariant map  
$$\mu_t: G \times^P (t + n(\p)) \to \g$$ defined by $\mu_t([g,t+x]) = Ad_g(t+x)$. 
The image $\mathrm{Im}(\mu_t)$ coincides with the closure $\bar{O}_t$ of an adjoint orbit $O_t$.  
By the property (2) of [Ibid, Proposition (1.4.14), (2)], $\omega_t$ is the pull-back of 
the Kostant-Kirillov form on $O_t$ by $\mu_t$. 
Denote by $(\widetilde{G\cdot r(\p) \times_{\h/W}\k(\p)})_{reg}$ the smooth part 
of $\widetilde{G\cdot r(\p) \times_{\h/W}\k(\p)}$. There is a crepant resolution    
$$ G \times^P r(\p) \to \widetilde{G\cdot r(\p) \times_{\h/W}\k(\p)}$$ which does not change the 
smooth locus. Then   
$\omega$ determines a relative symplectic 2-form $\bar{\omega}$ 
of $(\widetilde{G\cdot r(\p) \times_{\h/W}\k(\p)})_{reg} \to \k(\p)$.   
Note that $$\widetilde{G\cdot r(\p) \times_{\h/W} \k(\p)} \cong \widetilde{G\cdot r(\p)} \times_{\k(\p)/W'} \k(\p)$$ 
by Lemma (2.5). Since $W'$ acts on $\k(\p)$, it acts on $\widetilde{G\cdot r(\p)} \times_{\k(\p)/W'} \k(\p)$ in 
a natural manner.  
Let us check that $\bar{\omega}$ is $W'$-invariant. 
Take a general point $t \in \k(\p)$ (or more precisely, take $t \in \k(\p)^{reg}$ in the notation of 
[Na 4], (P1)). Then the fiber of the map 
$$\bar{f}: \widetilde{G\cdot r(\p)} \times_{\k(\p)/W'} \k(\p) \to \k(\p)$$ over $t$ is  
the semi-simple (adjoint) orbit $G \cdot t$ (cf. [Na 4], Lemma (1.1)). 
Then $\bar{\omega}_t$ coincides with the Kostant-Kirillov form on $G \cdot t$. 
If $w \in W'$, then $w(t) \in \k(\p)^{reg}$ and the fiber $\bar{f}^{-1}(w(t))$ is  
$G \cdot w(t)$. Moreover $\bar{\omega}_{w(t)}$ coincides with the Kostant-Kirillov form on 
$G \cdot w(t)$. Since $G\cdot t = G \cdot w(t)$, we see that $\bar{\omega}_t$ and $\bar{\omega}_{w(t)}$ 
coincides. This argument shows that $\bar{\omega}$ is $W'$-invariant on a Zariski open subset 
of $(\widetilde{G\cdot r(\p)} \times_{\k(\p)/W'} \k(\p))_{reg}$. Therefore $\bar{\omega}$ is $W'$-invariant. 
The relative 2-form $\bar{\omega}$ descends to a relative symplectic 2-form of $(\widetilde{G\cdot r(\p)})_{reg} 
\to \k(\p)/W'$, which determines a Poisson structure of $(\widetilde{G\cdot r(\p)})_{reg}$ 
over $\k(\p)/W'$. This Poisson structure uniquely extends to that of $\widetilde{G \cdot r(\p)}$. Q.E.D.  
\vspace{0.2cm}

{\bf Proposition (2.7)}. {\em The Poisson deformation $f: G \times^P r(\p) \to \k(\p)$ 
of $T^*(G/P)$ is universal  
at $0 \in \k(\p)$.} 
\vspace{0.2cm}

{\em Proof}. 
For $\lambda \in \k(\p)$, the vector space $\mathbf{C}\lambda \oplus n(\p)$ becomes a 
$P$-module by the adjoint action. Thus, one can define a vector bundle  
$G \times^P (\mathbf{C}\lambda \oplus n(\p))$ over $G/P$.  This vector bundle fits into the  
exact sequence
$$ 0 \to G \times^P n(\p) \to G \times^P (\mathbf{C}\lambda \oplus n(\p)) \to G/P \times \mathbf{C}\lambda 
\to 0.$$    
Since $G \times^P n(\p) \cong T^*(G/P)$, this extension gives an element $e(\lambda) \in H^1(G/P, \Omega^1_{G/P})$. 
Assume that $b_2(G/P) = n$ or equivalently that $\dim \k(\p) = n$.  
Then one can find $n$ maximal parabolic subgroups $Q_i$ $(1 \le i \le n)$ such 
that $P \subset Q_i$ and each projection map $\pi_i: G/P \to G/Q_i$ determines a extremal ray of the 
nef cone $\overline{\mathrm{Amp}}(G/P)$, which is a simplicial polyhedral cone of dimension $n$. 
Note that each $Q_i$ is the maximal parabolic subgroup associated to a vertex in the set of 
marked vertices corresponding to $P$ (cf. [Na 4], (P.3)).  
Let $\lambda_i \in \k(\q_i) \subset \k(\p)$ be a non-zero element. Note that, since $\k(\q_i)$ is one-dimensional, $\lambda_i$ 
is unique up to constant.  Then one has an exact sequence of vector bundles on $G/Q_i$: 
$$ 0 \to G \times^{Q_i} n(\q_i) \to G \times^{Q_i} (\mathbf{C}\lambda \oplus n(\q_i)) \to G/Q_i \times \mathbf{C}\lambda_i  
\to 0.$$ 
Let $f_i: G \times^{Q_i} (\mathbf{C}\lambda_i \oplus n(\q_i)) 
\to \k(\q_i)$ be the map defined by $f_i([g, t\lambda_i + x]) := t\lambda_i$ with 
$x \in n(\q_i)$. The fiber $f_i^{-1}(0) = T^*(G/Q_i)$ is not an affine variety, but 
$f_i^{^1}(\lambda_i)$ is an affine variety. Thus, the exact sequence above does not split, and 
its extension class $e(\lambda_i) \in H^1(G/Q_i, \Omega^1_{G/Q_i})$ is not zero. 
The exact sequence is pulled back to the exact sequence of vector bundles on $G/P$:  
$$ 0 \to G \times^P n(\q_i) \to G \times^P (\mathbf{C}\lambda \oplus n(q_i)) \to G/P \times \mathbf{C}\lambda_i \to 
0.$$  By the natural injection $G \times^P n(\q_i) \to G \times^P n(\p)$, one obtains an exact sequence  
$$ 0 \to G \times^P n(\p) \to G \times^P (\mathbf{C}\lambda_i \oplus n(\p)) \to G/P \times \mathbf{C}\lambda_i  
\to 0.$$ The extension class of this exact sequence is the image of $e(\lambda_i) \in H^1(G/Q_i, \Omega^1_{G/Q_i})$ 
by the map $H^1(G/Q_i, \Omega^1_{G/Q_i}) \to H^1(G/P, \Omega^1_{G/P})$.  
Note that $\{e(\lambda_i)\}$ is a basis of $H^1(G/P, \Omega^1_{G/P})$.  Therefore, the exact 
sequence 
$$ 0 \to G \times^P n(\p) \to G \times^P r(\p) \to G/P \times \k(\p) 
\to 0 $$ is the universal extension of  $G \times^P n(\p)$ by the trivial line bundle. 
Let $$ 0 \to \Omega^1_{G/P} \to \mathcal{E} \to \mathcal{O}_{G/P}^n \to 0$$ 
be the corresponding exact sequence of the sheaves.  
Let $p: T^*(G/P)(= G \times^P n(\p)) \to G/P$ be the canonical projection. 
Then we have a commutative diagram of exact sequences: 

\begin{equation}
\begin{CD}
0 @>>> p^*\Omega^1_{G/P} @>>> 
p^*\mathcal{E} @>>> p^*\mathcal{O}_{G/P}^{n} 
@>>> 0 \\ 
@. @VVV @VVV @VVV @. \\ 
0 @>>> \Theta_{T^*(G/P)} @>>> 
\Theta_{G \times^P r(\p)}\vert_{T^*(G/P)}
 @>>> 
N_{T^*(G/P)/G \times^P r(\p)} @>>> 0 
\end{CD}
\end{equation}
 
Identify $\k(\p)$ with its tangent space at $0$.   
The Kodaira-Spencer map $\theta_f$
of $f$ is given by the composite 
$$ \k(\p) \to H^0(T^*(G/P), N_{T^*(G/P)/G \times^P r(\p)}) 
\to H^1(T^*(G/P), \Theta_{T^*(G/P)}).$$ 
On the other hand, by the identification of $\k(\p)$ with 
$H^0(G/P, \mathcal{O}_{G/P}^n)$, one has a map 
$$ \k(\p) \cong H^0(G/P, \mathcal{O}_{G/P}^n) \to 
H^1(G/P, \Omega^1_{G/P}).$$ By the construction, the Kodaira-Spencer map 
is factored by this map: 
$$  \k(\p) \to H^1(G/P, \Omega^1_{G/P}) \to H^1(T^*(G/P), \Theta_{T^*(G/P)}). $$ 
The first map is an isomorphism by the definition of 
$\mathcal{E}$. The second map is an injection. In fact, let 
$S \subset T^*(G/P)$ be the zero section. Then $N_{S/T^*(G/P)} \cong 
\Omega^1_S$ and the composite $H^1(G/P, \Omega^1_{G/P}) 
\to H^1(T^*(G/P), \Theta_{T^*(G/P)}) \to H^1(S, \Omega^1_S)$ is an 
isomorphism. Therefore, the Kodaira-Spencer map 
$\theta_f$ is an injection. Since $f$ is a Poisson deformation 
of $T^*(G/P)$, the Kodaira-Spencer map $\theta_f$ is factored by 
the ``Poisson Kodaira-Spencer map"$\theta^P_f$: 
$$ \k(\p) \stackrel{\theta^P_f}\to H^2(T^*(G/P), \mathbf{C}) \to 
H^1(T^*(G/P), \Omega^1_{T^*(G/P)}).$$ Hence $\theta^P_f$ is also 
injective. Since $\dim \k(\p) = h^2(T^*(G/P), \mathbf{C}) = n$, 
$\theta^P_f$ is actually an isomorphism. Q.E.D.  
\vspace{0.2cm}

{\bf Theorem (2.8)}. {\em The Poisson deformation 
$\widetilde{G \cdot r(\p)} \to \k(\p)/W'$ is universal at $0 \in \k(\p)/W'$.} 
\vspace{0.2cm}

{\em Proof}. By Proposition (2.7), $G \times^P r(\p) \to \k(\p)$ is the 
universal Poisson deformation of $T^*(G/P)$ around $0 \in \k(\p)$. 
Since $\tilde{O}$ has a $\mathbf{C}^*$-action with positive weight, 
the universal Poisson deformation of $\tilde{O}$ is algebraized to a 
$\mathbf{C}^*$-equivariant map  
$\mathcal{Y} \to \mathbf{A}^d$ with $\mathcal{Y}_0 = \tilde{O}$. 
There is a $\mathbf{C}^*$-equivariant commutative diagram 
\begin{equation}
\begin{CD} 
G \times^P r(\p) @>>> \mathcal{Y}\\ 
@VVV @VVV \\ 
\k(\p) @>{\pi}>> \mathbf{A}^d. 
\end{CD}
\end{equation}  

If $t \in \k(\p)$ is general, then the induced map 
$G \times^P (t + n(\p)) \to \mathcal{Y}_{\pi(t)}$ is an isomorphism.    
By the main theorem, $\k(\p) \to \mathbf{A}^d$ is a finite Galois map. 
Denote by $H$ its Galois group. 
By Proposition (2.6), we have seen that $\widetilde{G\cdot r(\p)} \to \k(\p)/W'$ 
is a Poisson deformation of $\tilde{O}$. By the (formal) universality of 
$\mathcal{Y} \to \mathbf{A}^d$ at $0$, there is a formal map $\widehat{\k(\p)/W'} 
\to \hat{\mathbf A}^d$. Since 
$\widetilde{G\cdot r(\p)}$ has a $\mathbf{C}^*$-action and the Poisson deformation 
$\widetilde{G\cdot r(\p)} \to \k(\p)/W'$ is $\mathbf{C}^*$-equivariant, the 
formal map is also $\mathbf{C}^*$-equivariant. Then the formal map determines 
a map $\k(\p)/W' \to \mathbf{A}^d$ (cf. [Na 1], \S 7).  
The map $\k(\p) \to \mathbf{A}^d$ in the commutative diagram factorizes as  
$$\k(\p) \to \k(\p)/W' \to \mathbf{A}^d.$$ 
In fact, $$\mathrm{Spec}\Gamma(G \times^P r(\p), \mathcal{O}_{G \times^P r(\p)}) 
\to \k(\p)$$ is a Poisson deformation of $\tilde{O}$ and it coincides with 
the pull-back of $\mathcal{Y} \to \mathbf{A}^d$ by $\pi$. 
But $$\mathrm{Spec}\Gamma(G \times^P r(\p), \mathcal{O}_{G \times^P r(\p)}) 
= \widetilde{G\cdot r(\p)} \times_{\k(\p)/W'} \k(\p)$$ by Lemma (2.5). 
As a result, we have a commutative diagram of Poisson deformations of $\tilde{O}$:

\begin{equation}
\begin{CD} 
\mathrm{Spec}\Gamma(G \times^P r(\p), \mathcal{O}_{G \times^P r(\p)}) @>>> \widetilde{G\cdot r(\p)} 
@>>> \mathcal{Y} \\ 
@VVV @VVV @VVV\\ 
\k(\p) @>>> \k(\p)/W' @>>> \mathbf{A}^d. 
\end{CD}
\end{equation}  
   
We shall prove that $H = W'$. We only have to show that 
$H \subset W'$. An element $h \in H$ induces an automorphism of  
$\k(\p)$ over $\mathbf{A}^d$.  Denote also by $h$ this automorphism. 
Let us consider the automorphism of 
$\mathcal{Y} \times_{\mathbf{A}^d}\k(\p)$ defined by $id \times h$. 
It induces a birational automorphism $G \times^P r(\p) --\to G \times^P r(\p)$ 
and a commutative diagram 
$$G \times^P r(\p) --\to G \times^P r(\p)$$ 
$$\downarrow \hspace{1.5cm} \downarrow $$ 
$$\mathcal{Y} \times_{\mathbf{A}^d}\k(\p) \to \mathcal{Y} \times_{\mathbf{A}^d}\k(\p).$$  
Let $\mathcal{X}$ be the fiber product of the diagram 
$$\mathcal{Y} \times_{\mathbf{A}^d}\k(\p) \stackrel{id \times h}
\to \mathcal{Y} \times_{\mathbf{A}^d}\k(\p) 
\leftarrow G \times^P r(\p).$$ 
Then there is a birational map (over $\mathcal{Y} \times_{\mathbf{A}^d}\k(\p)$) 
$$ G \times^P r(\p) --\to \mathcal{X},$$ which is an isomorphism 
in codimension one.  Since $\mathcal{Y} \times_{\mathbf{A}^d}\k(\p)$ is affine, 
its normalization is isomorphic to $Y_{\k(\p)}$.   
Denote by $\mu_{\p}$ (resp. $\mu$) the birational map from 
$G \times^P r(\p)$ (resp. $\mathcal{X}$) to 
$Y_{\k(\p)}$. 
In $H^2(G \times^P r(\p), \mathbf{R})$, the movable cone 
$\mathrm{Mov}(\mu_0)$ 
coincides with $w(\mathrm{Mov}(\mu_{\p,0}))$ for some $w \in W'$. 
Let us consider the automorphism of 
$\mathcal{Y} \times_{\mathbf{A}^d}\k(\p)$ defined by 
$id \times h\cdot w^{-1}$. Let $\mathcal{X}'$ be the fiber product of the diagram 
$$\mathcal{Y} \times_{\mathbf{A}^d}\k(\p) \stackrel{id \times h\cdot w^{-1}}
\to \mathcal{Y} \times_{\mathbf{A}^d}\k(\p) 
\leftarrow G \times^P r(\p).$$ As above, there is a birational map (over 
$\mathcal{Y} \times_{\mathbf{A}^d}\k(\p)$; hence over $Y_{\k(\p)}$) 
$$ G \times^P r(\p) --\to \mathcal{X}'.$$  
Denote by $\mu'$ the birational map 
from $\mathcal{X}'$ to $Y_{\k(\p)}$. 
Then, by the construction, $\mathrm{Mov}(\mu'_0)$ 
coincides with $\mathrm{Mov}(\mu_{\p,0})$. 
It follows that $G \times^P r(\p)$ and $\mathcal{X}'$ are connected by 
a sequence of flops of the 1-st kind: 
$$ G \times^P r(\p) --\to G \times^{P_1}r(\p_1) --\to ... --\to 
G \times^{P_k} r(\p_k) = \mathcal{X}'.$$
Let $E \subset T^*(G/P)$ be an exceptional divisor 
of the Springer map $s: T^*(G/P) \to \tilde{O}$. At a general point of $E$, 
the flop $G \times^P r(\p) --\to G \times^{P_1} r(\p_1)$ is an isomorphism. 
Let $E_1$ be the proper transform of $E$ by this flop. At a general point of $E_1$, 
the next flop $G \times^{P_1} r(\p_1) --\to G \times^{P_2} r(\p_2)$ is also an 
isomorphism. Similar things happen for all flops of the 1-st kind.  
If $h\cdot w^{-1} \ne 1$, then, by the main theorem, the indeterminacy locus 
of the birational map $G \times^P r(\p) --\to \mathcal{X}'$ must contain 
at least one $s$-exceptional divisor $E$ (cf. Remark below (1.1)). Therefore, $h\cdot w^{-1} = 1$. 
Q.E.D.  
\vspace{0.2cm}

{\bf Example (2.9)}. The following are standard examples due to Slodowy [Slo]. 
Let $\g$ be of type $B_m$, $C_m$, $F_4$ or $G_2$. Let $B$ be a Borel 
subgroup of $G$. Then the Springer map $s: T^*(G/B) \to \g$ gives a 
crepant resolution of the nilpotent cone $N$; namely, $\mathrm{Im}(s) = N$. 
Note that $N$ is the closure of the regular nilpotent orbit of $\g$. 
In this case, $\Sigma - \Sigma_0$ has only one connected component, say $B$. 
Then the surface $S_B$ has a singularity of type $A_{2m-1}$, $D_{m+1}$, $E_6$ or 
$D_4$ according as $\g$ is of type $B_m$, $C_m$, $F_4$ or $G_2$. 
Note that $W'$ is the Weyl group $W(G)$ of $G$. By Theorem 2.8 and Theorem 1.1, 
we see that $W_B = W(G)$. This means that $\mathrm{Exc}(s)$ is a divisor of 
$T^*(G/B)$ with exactly  $m$ (resp. $4$, $2$) irreducible components when $\g$ is of type $B_m$ or 
$C_m$ (resp. $F_4$, $G_2$).     
\vspace{0.2cm}

\begin{center}
{\bf References} 
\end{center}

[Be] Beauville, A.: Vari\'{e}t\'{e}s K\"{a}hleriennes dont la premi\`{e}re classes de Chern est nulle, 
J. Diff. Geom. {\bf 18} (1983) 755-782 

[Ca] Carter, R.: Simple groups of Lie type, John Wiley and Sons 1972. 

[Ca 2] Carter, R.: Lie algebras of finite and affine type, Cambridge studies in advanced 
mathematics {\bf 96},  2005 

[C-G] Chriss, N., Ginzburg, V.: Representation theory and complex geometry, 
Birkhauser, Boston, MA, 1997  

[G-K] Ginzburg, V., Kaledin, D.: Poisson deformations of symplectic quotient 
singularities, Adv. Math. {\bf 186} (2004), 1-57

[Ho] Howlett, R.B.: Normalizers of parabolic subgroups of reflection 
groups, J. London Math. Soc. {\bf 21} (1980), 62-80 

[Ka] Kaledin, D.: Symplectic singularities from the Poisson point of view, 
J. Reine Angew. Math, {\bf 600} (2006) 135-156 

[Ko-Mo] Koll\'{a}r, J., Mori,S.: Classification of the three dimensional 
flips, J. Amer. Math. Soc. {\bf 5} (3) (1992), 533-703

[Lo] Looijenga, E.: Isolated singular points on complete 
intersections, London Math. Soc. Lecture Note Series {\bf 77}, 
Cambridge University Press (1984)

[Ma] Markman, E.: Modular Galois covers associated to symplectic 
resolutions of singularities, arxiv: 0807.3502

[Na] Namikawa, Y.: Deformation theory of singular symplectic n-folds, 
Math. Ann. {\bf 319} (2001), 597-623 

[Na 1] Namikawa, Y.: Flops and Poisson deformations of symplectic 
varieties, Publ. RIMS, Kyoto Univ. {\bf 44} (2008), 259-314 

[Na 2] Namikawa, Y.: Poisson deformations of affine symplectic varieties, 
math.AG/0609741 

[Na 3] Namikawa, Y.: On deformations of {\bf Q}-factorial symplectic varieties, 
J. Reine Angew. Math. {\bf 599} (2006), 97-110 

[Na 4]Namikawa, Y.: Birational geometry and deformations of nilpotent 
orbits, Duke Math. J. {\bf 143} (2008), 375-405 

[Na 5] Namikawa, Y.: ``Birational geometry of symplectic resolutions of 
nilpotent orbits'' in {\em Moduli spaces and 
Arithmetic geometry (Kyoto 2004)}, Advvanced Studies in Pure Mathematics {\bf 45} (2006) 
75-116

[Na 6] Namikawa, Y.: On deformations 
of {\bf Q}-factorial symplectic varieties, J. Reine Angew. Math. 
{\bf 599} (2006) 97-110 

[Slo] Slodowy, P.: Simple singularities and simple algebraic groups, 
Lecture Notes in Math. {\bf 815}, Springer-Verlag (1980) 

[Slo 2] Slodowy, P.: Four lectures on simple groups and singularities, 
Comm. of Math. Inst. Rijksuniversiteit Utrecht 1980 

[Ya] Yamada, H.: Lie group theoretical construction of period mapping, 
Math. Z. {\bf 220} (1995), 231-255 

\vspace{0.3cm}

\begin{center}
Department of Mathematics, Faculty of Science, Kyoto University, Japan 

namikawa@math.kyoto-u.ac.jp
\end{center}

\end{document}